\documentclass[10pt]{amsart}

\usepackage{amssymb}
\usepackage{amscd}
\usepackage{upref}
\usepackage{tikz} 
\usetikzlibrary{decorations.pathmorphing}
\usepackage[all]{xy}

\theoremstyle{plain}
\newtheorem{thm}{Theorem}[section]
\newtheorem{cor}[thm]{Corollary}
\newtheorem{lem}[thm]{Lemma}
\newtheorem{prop}[thm]{Proposition}

\newtheorem*{Theorem A}{Theorem A}
\newtheorem*{Theorem B}{Theorem B}
\newtheorem*{Theorem C}{Theorem C}
\newtheorem*{Theorem D}{Theorem D}

\theoremstyle{definition}
\newtheorem{defn}[thm]{Definition}
\newtheorem{example}[thm]{Example}

\theoremstyle{remark}

\newtheorem{rem}[thm]{Remark}

\numberwithin{equation}{section}

\begin{document}
\title[Open-Closed Disc String Operations and Saddle Interactions]{Open-Closed TQFT String Operations for Disc Cobordisms, Simultaneous Saddle Interactions, and Constant Homology Classes}
\author{Hirotaka Tamanoi}
\address[]
{Department of Mathematics, University of California Santa Cruz \newline
\indent Santa Cruz, CA 95064}
\email[]{tamanoi@math.ucsc.edu}
\date{}
\begin{abstract}
In \cite{T5}, we showed that most of the open-closed topological quantum field theory (TQFT) string operations vanish including all the higher genus TQFT operations, and we described a small list of genus zero open-closed TQFT string operations which can be nontrivial. In this paper, we consider open-closed string operations associated to open-closed cobordisms homeomorphic to discs. These operations constitute the main part of genus zero string operations, and they include the saddle string operation of two open strings interacting at their internal points. We show not only that disc string operations are independent of their half-pair-of-pants decompositions but also that these disc string operations can be computed by \emph{simultaneous} saddle interactions of incoming open strings at the same point, and they take values in homology classes of constant open strings on some closed orientable submanifolds, which we will precisely determine. We will also discuss a role played by fundamental constant homology classes in open-closed string topology. Our main tools are saddle interaction diagrams and their deformations. 
\end{abstract}
\maketitle

\tableofcontents

\section{Introduction}

Let $M$ be an oriented closed connected smooth manifold of dimension $d$ and let $LM$ be its free loop space of continuous loops in $M$. Chas and Sullivan \cite{CS} discovered the Batalin-Vilkovisky algebra structure in the homology of the free loop space $\mathbb H_*(LM)=H_{*+d}(LM)$, and started string topology. Subsequently, string operations in closed string topology were constructed using homotopy theoretic methods in \cite{CJ} and also in \cite{CG}. In this paper, we discuss string operations in open-closed string topology, continuing our discussion in \cite{T5}. 

Let $\Sigma$ be an orientable 2-dimensional open-closed cobordism surface in open-closed string topology whose free boundary components are labeled by oriented closed submanifolds $I,J,K,L,\dots$ of $M$. It is assumed that $\Sigma$ has at least one outgoing closed or open string (positive boundary condition). Let $P_{IJ}$ be the space of open strings consisting of continuous maps $\gamma:[0,1] \to M$ starting at points in $I$ and ending at points in $J$. Suppose $\Sigma$ has $p$ incoming and $q$ outgoing closed strings, and suppose $(I,J)$'s are labels of incoming open strings and $(K,L)$'s are labels of outgoing open strings. Then open-closed string topology associates an operator $\mu_{\Sigma}$, the TQFT string operation associated to the open-closed cobordism $\Sigma$, of the following form (modulo K\"unneth Theorem),  
\begin{equation*}
\mu_{\Sigma}: H_*(LM)^{\otimes p}\otimes{\textstyle\bigotimes_{(I,J)}} H_*(P_{IJ}) 
\longrightarrow H_*(LM)^{\otimes q}\otimes\textstyle{\bigotimes_{(K,L)}}
H_*(P_{KL}),
\end{equation*}
determined up to sign (\cite{Ra}, \cite{Su}). See also \cite{Harrelson} for genus $0$ open-closed TQFT, and \cite{LP} for a categorical approach to open-closed TQFT of arbitrary genera. In \cite{Go2}, Godin constructed higher string operations associated to homology classes of the mapping class groups of open-closed cobordism surfaces, when only one label is used in open-closed string topology. The TQFT operations in this paper correspond to degree zero homology classes of the mapping class groups, and they depend only on the homeomorphism types of open-closed cobordism surfaces with labeled free boundaries. For more information on higher string operations, see \cite{T7} where we showed that stable higher string operations are all trivial. We often omit the phrase TQFT and simply refer to them as string operations. 

In \cite{T3} and \cite{T5}, we showed that most of the open-closed TQFT string operations, including all the higher genus ones, vanish, except for string operations associated to the following open-closed cobordisms $\Sigma$ of genus zero: 
\begin{enumerate} 
\item $\Sigma$ has genus $0$, one window, with no open strings, and with one outgoing closed string. 
\item $\Sigma$ has genus $0$, with no windows, no outgoing open strings, and with $1$ or $2$ outgoing closed strings. 
\item $\Sigma$ has genus $0$, no windows, no outgoing closed strings, and all outgoing open strings are along exactly one boundary component of $\Sigma$ which may also contain incoming open strings. 
\end{enumerate} 
In Theorem B and in Figures 1 to 5 of \cite{T5}, we listed five types (I) to (V) of open-closed cobordisms with possibly nontrivial string operations. The above case (1) corresponds to (I) in Theorem B, the case (2) corresponds to (II) and (III) in Theorem B, and the case (3) corresponds to (IV) and (V) of Theorem B in that paper. 

The associated string operations in the cases (1) and (2) can be immediately described, provided that we explicitly know closed string products (loop products) and open string products for $M$. For the description of open window operation and closed window operation, see  Proposition 3.2 and Proposition 3.6 in \cite{T5}. With these informations, computing string operations in case (3) reduces to understanding the TQFT string operations associated to open-closed cobordisms homeomorphic to discs, which is the purpose of this paper. 

Of course, we can compute the string operation associated to a disc cobordism $\Sigma$ by decomposing $\Sigma$ into half-pair-of-pants, and then composing open string products and coproducts according to the  decomposition of $\Sigma$. Here, a half-pair-of-pants is an open cobordism of the form 
\begin{tikzpicture}
\draw[scale=0.1, rotate around={180:(0,0.5)}] (0,0) -- ++(1,0) arc (90:45:1.5) arc (225:270:1.5) --  ++(1,0)    -- ++(0,1) -- ++(-1,0) arc (270:90:0.5) -- ++(1,0)-- ++(0,1)
 -- ++(-1,0) arc (90:135:1.5) arc (315:270: 1.5) -- ++(-1,0) -- cycle;
\end{tikzpicture}
or
\begin{tikzpicture}
\draw[scale=0.1] (0,0) -- ++(1,0) arc (90:45:1.5) arc (225:270:1.5) --  ++(1,0)    -- ++(0,1) -- ++(-1,0) arc (270:90:0.5) -- ++(1,0)-- ++(0,1)
 -- ++(-1,0) arc (90:135:1.5) arc (315:270: 1.5) -- ++(-1,0) -- cycle;
\end{tikzpicture}
whose corresponding string operations are the open string product and the coproduct. 
However, precise understanding of the open string product and coproduct is not so easy. Furthermore, if we compose open string products and coproducts according to a pair-of-pants decomposition of a disc cobordism $\Sigma$, we can easily imagine that the resulting string operation can be very complicated. 

In this paper, we aim at describing qualitative properties of disc string operations. See Theorem B and Theorem C below. However, one case in which an explicit computation is possible is when $\Sigma$ has only outgoing open strings and no incoming open strings. Let $k$ be the coefficient field of homology. 

\begin{figure}
\begin{center}
\begin{tikzpicture}
\draw (0,0) circle (1cm);
\draw (35:0.97cm) -- (35:1.03cm) arc (35:55:1.03cm) 
      -- (55:0.97cm) arc (55:35:0.97cm);
\draw (125:0.97cm) -- (125:1.03cm) arc (125:145:1.03cm) 
      -- (145:0.97cm) arc (145:125:0.97cm);
\draw (350:0.97cm) -- (350:1.03cm) arc (350:360:1.03cm) 
 arc (0:10:1.03cm) -- (10:0.97cm) arc (10:0:0.97cm)
 arc (360:350:0.97cm);
\draw (170:0.97cm) -- (170:1.03cm) arc (170:190:1.03cm) 
      -- (190:0.97cm) arc (190:170:0.97cm);
\draw (215:0.97cm) -- (215:1.03cm) arc (215:235:1.03cm) 
      -- (235:0.97cm) arc (235:215:0.97cm);
\fill (35:1cm) circle (1pt);
\fill (55:1cm) circle (1pt);
\fill (125:1cm) circle (1pt);
\fill (145:1cm) circle (1pt);
\fill (350:1cm) circle (1pt);
\fill (10:1cm) circle (1pt);
\fill (170:1cm) circle (1pt);
\fill (190:1cm) circle (1pt);
\fill (215:1cm) circle (1pt);
\fill (235:1cm) circle (1pt);

\path (0:1.1cm) node[right] {$\eta_r$};
\path (45:1cm) node[above right] {$\eta_{r+1}$};
\path (135:1cm) node[above left] {$\eta_1$};
\path (180:1.1cm) node[left] {$\eta_2$};
\path (225:1cm) node[below left] {$\eta_3$};

\draw[->,>=stealth,  thick] (2:1cm) -- (4:1cm);
\draw[->,>=stealth,  thick] (227:1cm) -- (229:1cm);
\draw[->,>=stealth,  thick] (182:1cm) -- (184:1cm);
\draw[->,>=stealth,  thick] (137:1cm) -- (139:1cm);
\draw[->,>=stealth,  thick] (47:1cm) -- (49:1cm);
\path (90:1cm) node[above] {$I$};
\path (25:1cm) node[right] {$K_r$};
\path (150:1cm) node[left] {$K_1$};
\path (210:1cm) node[left] {$K_2$};
\draw[dashed] (250:1.15cm) arc (250:340:1.15cm);
\draw node[text width=10.5cm] at (0,-2.8) 
{\textsc{Figure 1.} A disc cobordism $\Sigma$ with no incoming open strings and $r+1$ outgoing open strings $\eta_1,\eta_2,\dots,\eta_{r+1}$ having arc labels $K_1,K_2,\dots,K_r$ between open strings. The associated string operation has values in homology classes of constant open strings. See Theorem A.  };
\end{tikzpicture}
\end{center}
\end{figure}

\begin{Theorem A} Let $\Sigma$ is an open-closed cobordism homeomorphic to a disc with $r+1$ outgoing open strings and no incoming open strings with $r\ge0$. Suppose that its free boundaries are labeled by closed oriented mutually transversal submanifolds $I, K_1, \dotsc, K_r$, where $I$ is the outer most label of $\Sigma$ \textup{(}See Figure 1\textup{)}. Then the associated string operation $\mu_{\Sigma}$ is given by 
\begin{align*}
\mu_{\Sigma}:&k \longrightarrow H_*(P_{IK_1})\otimes H_*(P_{K_1K_2})\otimes\dotsb\otimes H_*(P_{K_rI}), \\
\mu_{\Sigma}(1)&=\pm (s_*\otimes\dotsb\otimes s_*)\circ (\phi_{r+1})_*
\bigl([I\cap K_1\cap\dotsb\cap K_r]\bigr),
\end{align*}
where $\phi_{r+1}$ is the iterated diagonal map into $(r+1)$-fold  Cartesian product of $I\cap K_1\cap\dotsb\cap K_r$, and $s$ the inclusion map from $I\cap K_1\cap\dotsb\cap K_r$ into $P_{K_iK_{i+1}}$ for $0\le i\le r$ with $K_0=K_{r+1}=I$. 
\end{Theorem A}

Note that the image $\mu_{\Sigma}(1)$ is a homology class of constant paths. This constancy of resulting homology classes turns out to be always the case for (iterated) open string coproducts and general string operations with at least two outgoing strings. Moreover, open string products with basic constant homology classes $[I\cap J]\in H_*(P_{IJ})$ always result in constant homology classes, and they are closely related to the first and the last terms in the open string coproduct. 

More precisely, let $\varphi_J: H_*(P_{IK}) \rightarrow H_*(P_{IJ})\otimes H_*(P_{JK})$ be the open string coproduct map in which open strings are split along the submanifold $J$, and let $\mu_J: H_*(P_{IJ})\otimes H_*(P_{JK}) \rightarrow H_*(P_{IK})$ be the open string product along $J$. 

\begin{Theorem B} \textup{(1)} The open string coproduct map $\varphi_J$ factors through homology groups of constant open strings, as follows. 
\begin{equation*}
\xymatrix{
H_*(P_{IK}) \ar[r]^{\varphi_J\ \ \ \ \ \ \ } \ar[dr] &
H_*(P_{IJ})\otimes H_*(P_{JK}) \\
& H_*(I\cap J)\otimes H_*(J\cap K) \ar[u]_{s_*\otimes s_*}
}
\end{equation*}
Here $s_*$ is induced from an appropriate inclusion map. 

\textup{(2)}  For $a\in H_*(P_{IK})$, the first and the last terms in the open string coproduct $\varphi_J(a)$ is given by open string products with the basic constant homology classes $[J\cap K]$ and $[I\cap J]$ up to sign\textup{:}
\begin{equation*}
\varphi_J(a)=\mu_K(a\otimes[J\cap K])\otimes(\pm1)+\cdots+(\pm1)\otimes \mu_I([I\cap J]\otimes a),
\end{equation*}
where $\mu_K(a\otimes[J\cap K])\in H_*(P_{IJ})$ and $\mu_I([I\cap J]\otimes a)\in H_*(P_{JK})$ are certain homology classes of constant open strings. 
\end{Theorem B}

For a precise statement on the nature of the open string product with basic constant homology classes $[I\cap J]\in H_*(P_{IJ})$, see part (2) in Theorem \ref{image of open coproducts}.

\begin{figure}
\begin{center}
\begin{tikzpicture}
\draw (0,0) circle (1cm);
\draw[ultra thick] (0:1cm) arc (0:10:1cm);
\draw[ultra thick] (350:1cm) arc (350:360:1cm);
\draw[ultra thick] (80:1cm) arc (80:100:1cm);
\draw[ultra thick] (260:1cm) arc (260:280:1cm);
\draw[ultra thick] (170:1cm) arc (170:190:1cm);
\draw (35:0.97cm) -- (35:1.03cm) arc (35:55:1.03cm) 
      -- (55:0.97cm) arc (55:35:0.97cm);
\draw (125:0.97cm) -- (125:1.03cm) arc (125:145:1.03cm) 
      -- (145:0.97cm) arc (145:125:0.97cm);
\draw (305:0.97cm) -- (305:1.03cm) arc (305:325:1.03cm) 
      -- (325:0.97cm) arc (325:305:0.97cm);
\fill (10:1cm) circle (1pt);      
\fill (350:1cm) circle (1pt);      
\fill (80:1cm) circle (1pt);      
\fill (100:1cm) circle (1pt);      
\fill (260:1cm) circle (1pt);      
\fill (280:1cm) circle (1pt);      
\fill (170:1cm) circle (1pt);      
\fill (190:1cm) circle (1pt);      
\fill (35:1cm) circle (1pt);      
\fill (55:1cm) circle (1pt);      
\fill (125:1cm) circle (1pt);      
\fill (145:1cm) circle (1pt);      
\fill (305:1cm) circle (1pt);      
\fill (325:1cm) circle (1pt);      

\draw[->,>=stealth,  thick] (368:1cm) -- (356:1cm);
\draw[->,>=stealth,  thick] (88:1cm) -- (86:1cm);
\draw[->,>=stealth,  thick] (178:1cm) -- (176:1cm);
\draw[->,>=stealth,  thick] (268:1cm) -- (266:1cm);
\draw[->,>=stealth,  thick] (47:1cm) -- (49:1cm);
\draw[->,>=stealth,  thick] (317:1cm) -- (319:1cm);
\draw[->,>=stealth,  thick] (137:1cm) -- (139:1cm);
\draw[dashed] (200:1.2cm) arc (200:250:1.2cm);
\path (90:1cm) node[above] {$\gamma_1$};
\path (45:1cm) node[above right] {$\eta_1$};
\path (0:1cm) node[right] {$\gamma_2$};
\path (315:1cm) node[below right] {$\eta_2$};
\path (270:1cm) node[below] {$\gamma_3$};
\path (135:1cm) node[above left] {$\eta_m$};
\path (180:1cm) node[left] {$\gamma_m$};
\path (107:1cm) node[above left] {$K_1$};
\path (75:1cm) node[above right] {$K_2$};
\path (22:1cm) node[right] {$K_3$};
\path (330:1cm) node[right] {$K_4$};
\path (285:1cm) node[below right] {$K_5$};
\path (150:1cm) node[left] {$K_{2m}$};
\draw node[text width=11cm] at (0,-2.5) 
{\textsc{Figure 2.} A disc cobordism with alternating $m$ incoming open strings $\gamma_1,\gamma_2,\dots,\gamma_m$ and $m$ outgoing open strings $\eta_1,\eta_2,\dots,\eta_m$, with free arc labels $K_1,K_2,\dots, K_{2m}$ between open strings. The corresponding string operation is described in Theorem C. };
\end{tikzpicture}
\end{center}
\end{figure}

Although the open string product can be highly nontrivial, the above result shows that the open string coproduct behaves in a much simpler way. This is similar to the situation in closed string topology. In \cite{T3} we showed that the loop coproduct in $H_*(LM)$ behaves in a dramatically simpler way compared with the loop product: the coproduct is nontrivial only on $H_d(LM)$ with $d=\dim M$ and then it has values in constant homology classes in $H_0(LM)\otimes H_0(LM)$.  

Since the method of the proof of Theorem B, stated as Theorem 4.2 in \S 4, is homological, we also give geometric explanation of the above result of constant homology classes in terms of transverse intersection of cycles in the spirit of \cite{CS}. This explanation should be more intuitive and illuminating. See the discussion following the proof of Theorem 4.2. There, we also point out a role played by homology classes of constant closed strings in the $\mathbb H_*(LM)$-module structure on $H_*(P_{JK})$. See \eqref{LM-module structure}. 

We can state similar results for general iterated coproducts. See
Theorem \ref{iterated coproducts} in \S 4.2. The basic idea of the proof is to deform sewing diagrams (\S2) for iterated open string coproducts to degenerate configurations which manifest essential qualitative features of these operations. 

Now let $\Sigma$ be a general open-closed cobordism homeomorphic to a disc. We examine the associated string operation $\mu_{\Sigma}$. By use of open string products and open string coproducts, we are reduced to the case in which incoming open strings and outgoing open strings alternate along the boundary of $\Sigma$. To be more precise, let $K_1,K_2,\dots, K_{2m}$ be labels for free arcs along the boundary of $\Sigma$ separated by open strings. For convenience, let $P_{i,j}=P_{K_i,K_j}$. For $1\le k\le m$ suppose that $P_{2k-1,2k}$ and $P_{2k+1,2k}$ represent the configuration spaces of $k$-th incoming and $k$-th outgoing open strings. See Figure 2. Then the string operation $\mu_{\Sigma}$ is of the form 
\begin{equation}\label{string operation}
\mu_{\Sigma}: \bigotimes_{k=1}^mH_*(P_{2k-1,2k}) \longrightarrow \bigotimes_{k=1}^m H_*(P_{2k+1,2k}).
\end{equation}

We show that elements in the image of $\mu_{\Sigma}$ are not only homology classes of constant open strings, but also they are constant homology classes contained in some specific orientable submanifolds which depend only on $\Sigma$ and labeling submanifolds $K_1,K_2,\dots,K_{2m}$, and we describe them explicitly in Theorem C. This is a consequence of various possible interactions of the $m$ incoming open strings $\{\gamma_k\}_{1\le k\le m}$  at their internal points, not just at their end points. If we visualize two open strings moving in the manifold $M$, it is far more likely that they meet at their internal points, rather than at their end points, which occur as special cases. When two incoming open strings interact at their internal points for cutting and rejoining and become two outgoing open strings, they sweep out a surface which looks like a saddle. In Proposition 3.5 of \cite {T5}, this string operation is introduced and is called a saddle operation, and we showed that double saddle operations are always trivial. A saddle operation can be shown to be a composition of an open string coproduct followed by an open string product. For details, see the discussion of type (1) deformation of saddle interaction diagrams in \S 3.2. 

We introduce sewing diagrams in \S 2 to encode sequences of open string products and coproducts resulting from half-pair-of-pants decompositions of open-closed cobordisms. Note that in the open string product, two open strings interact at their end points. To handle general open string interactions at internal points, we introduce saddle interaction diagrams in \S 3 generalizing two-string saddle operations. A generic saddle interaction diagram describes a family of open strings interacting as pairs at their internal points. By deforming this configuration continuously, we can consider various \emph{simultaneous} saddle interaction processes of open strings in which all incoming open strings meet and interact at the same point, and the location of this interaction point on open strings is parametrized by $\vec {t}\in[0,1]^m$. See Figure 3 and diagram \eqref{simultaneous interaction} in \S5. It turns out that under a mild condition, the associated string operations remain the same under this continuous deformation of saddle interaction diagrams. Considering deformation into various simultaneous saddle interaction diagrams, we can reveal intrinsic qualitative topological properties of disc string operations which we describe in Theorem C below.

\begin{figure}
\begin{center}
\begin{tikzpicture}
\draw (-5.5,0.5) -- (-5.5,-0.5);
\draw (-4.5,0.5) -- (-4.5,-0.5);
\draw (-4,0.5) -- (-4,-0.5);
\draw (-3,0.5) -- (-3,-0.5);
\draw[dashed] (-5.3,0) -- (-4.7,0);
\draw[dashed] (-3.8,0) -- (-3.2,0);
\fill (-5.5,0.5) circle (2pt);
\fill (-5.5,-0.5) circle (2pt);
\fill (-4.5,0.5)circle (2pt);
\fill (-4.5,-0.5) circle (2pt);
\fill (-4,0.5) circle (2pt);
\fill (-4,-0.5) circle (2pt);
\fill (-3,0.5) circle (2pt);
\fill (-3,-0.5) circle (2pt);
\fill (-5.5,0) circle (1pt);
\path (-5.5,0.1) node[left] {$t_1$};
\fill (-4.5,0.2) circle (1pt);
\path (-4.5,0.2) node[right] {$t_k$};
\fill (-4,0.1) circle (1pt);
\path (-4,0.2) node[right] {$t_{k+1}$};
\fill (-3,0.1) circle (1pt);
\path (-3,0.1) node[right] {$t_m$};
\draw[->,>=stealth] (-5.5, -0.2) -- ++(0,-0.01);
\draw[->,>=stealth]  (-4.5,-0.2) -- ++(0,-0.01);
\draw[->,>=stealth]  (-4,-0.2) -- ++(0,-0.01);
\draw[->,>=stealth]  (-3,-0.2) -- ++(0,-0.01);
\path (-5.5,0.6) node[above] {$\gamma_1$};
\path (-4.5,0.6) node[above] {$\gamma_k$};
\path (-3.9,0.6) node[above] {$\gamma_{k+1}$};
\path (-3,0.6) node[above] {$\gamma_m$};

\draw[->,>=stealth] (-1.7,0) -- (-2.3,0);
\path (-2,0) node[above] {$\iota_{\vec{t}}$};

\draw (90:1cm) -- (0,0) -- (270:1cm) 
 (30:0.6cm) -- (0,0) -- (210:1.4cm) 
 (0:0.8cm) -- (0,0) -- (180:1.2cm) 
 (120:0.8cm) -- (0,0) -- (300: 1.2cm);
\draw[dashed] (40:0.8cm) arc (40:80:0.8cm);
\draw[dashed] (130:0.8cm) arc (130:170:0.8cm);
\draw[dashed] (220:0.8cm) arc (220:260:0.8cm);
\draw[dashed] (310:0.8cm) arc (310:350:0.8cm);
\fill (90:1cm) circle (2pt);
\fill (270:1cm) circle (2pt);
\fill (30:0.6cm) circle (2pt);
\fill (210:1.4cm) circle (2pt);
\fill(0:0.8cm) circle (2pt);
\fill (180:1.2cm) circle (2pt);
\fill(120:0.8cm) circle (2pt);
\fill (300: 1.2cm) circle (2pt);
\path (0:0.8cm) node[right] {$\gamma_{k+1}$};
\path (30:0.6cm) node[right] {$\gamma_k$};
\path (90:1cm) node[above] {$\gamma_1$};
\path (120:0.8cm) node[above left] {$\gamma_m$};
\draw[->,>=stealth] (0:0.4cm) -- (0:0.39cm);
\draw[->,>=stealth]  (30:0.3cm) -- (30:0.29cm);
\draw[->,>=stealth]  (90:0.5cm) -- (90:0.49cm);
\draw[->,>=stealth]  (120:0.4cm) -- (120:0.39cm);

\draw[->,>=stealth] (1.7,0) -- (2.3,0);
\path (2,0) node[above] {$j_{\vec{t}}$};

\draw (3,0.5) -- (3,-0.5);
\fill (3,0.5) circle (2pt);
\fill (3,-0.5) circle (2pt);
\fill (3,0) circle (1pt);
\draw[->, >=stealth] (3,-0.2) -- ++(0,-0.01);
\path (3,0.6) node[above] {$\eta_1$};
\draw[dashed] (3.2,0) -- (3.8, 0);
\draw (4,0.5) -- (4,-0.5);
\fill (4,0.5) circle (2pt);
\fill (4,-0.5) circle (2pt);
\path (4,0.6) node[above] {$\eta_k$};
\fill (4,0.1) circle (1pt);
\draw[->,>=stealth] (4,-0.2) -- ++(0,-0.01);
\path (4,0.3) node[right] {$\gamma_{k+1}|_{[0,t_{k+1}]}$};
\path (4,-0.3) node[right] {$\gamma_k|_{[t_k,1]}$};
\draw[dashed] (4.2, 0) -- (5.9,0);
\draw (6.2,0.5) -- (6.2,-0.5);
\fill (6.2,0.5) circle (2pt);
\path (6.2,0.6) node[above] {$\eta_m$};
\fill (6.2,-0.5) circle (2pt);
\fill (6.2,0) circle (1pt);
\draw[->, >=stealth] (6.2,-0.2) -- ++(0,-0.01);

\draw node[text width= 11cm] at (0.5,-3)
{\textsc{Figure 3.} A simultaneous saddle interaction of $m$ incoming open strings $\gamma_1,\dots,\gamma_m$ at the same point with parameter $\vec{t}=(t_1,\dots,t_m)\in [0,1]^m$ producing $m$ outgoing open strings $\eta_1,\dots,\eta_k,\dots,\eta_m$, where the first part of $\eta_k$ comes from the beginning portion of $\gamma_{k+1}$ and the second part comes from the end portion of $\gamma_k$. This interaction gives rise to a string operation $\mu_{\Sigma}(\vec{t})$ defined in \eqref{parametrized string operation}. See \S5 for further explanations. };
\end{tikzpicture}
\end{center}
\end{figure}

We note that the open-closed cobordisms corresponding to saddle interaction diagrams are disc cobordisms described in  Figure 2. 

As a special case of simultaneous saddle interactions, suppose incoming open strings $\{\gamma_k\}_{1\le k\le m}$ interact simultaneously at their end points (head or tail vertices carrying labels) rather than at their internal points, corresponding to corner points of the space $[0,1]^m$ parametrizing simultaneous saddle interactions. Here we treat open string interactions at end points not as open string products and coproducts, but as saddle interactions. For a discussion on the difference between sewing diagrams, which describe open string products and coproducts as special cases, and saddle interaction diagrams, see Remark \ref{difference between diagrams}. Suppose the interacting end point on the $k$-th incoming open string $\gamma_k\in P_{2k-1,2k}$ carries a label $K_{2k-1+\varepsilon_k}$ with $\varepsilon_k=0,1$ for $1\le k\le m$. The diagram for this simultaneous interaction is given in \eqref{e-interaction}. For a sequence $\varepsilon=(\varepsilon_1,\varepsilon_2,\dots,\varepsilon_m)$ with $\varepsilon=0,1$, let $K_{\varepsilon}=\bigcap_{k=1}^mK_{2k-1+\varepsilon_k}$. Figure 4 describes such an interaction when $\varepsilon=(1,\dots,1,0,\dots,0)$ with the last 1 at the $k$th position. In this saddle interaction, the simultaneous saddle interaction point on the manifold $M$ lies inside the orientable submanifold $K_{\varepsilon}$, and one of the outgoing open strings is a constant open string if $\varepsilon\not=(0,\dots,0), (1,\dots,1)$. 

\begin{Theorem C} The image of the string operation $\mu_{\Sigma}$ in \eqref{string operation} associated to a disc cobordism $\Sigma$ with alternating $m$ incoming and $m$  outgoing open strings with free arc labels $K_1,K_2,\dots, K_{2m}$ is contained in the following subgroup of constant homology classes\textup{:} 
\begin{equation*}
\textup{Im}\,\mu_{\Sigma}\subset \bigotimes_{k=1}^mS_{2k+1,2k}\subset \bigotimes_{k=1}^m H_*(P_{2k+1,2k}). 
\end{equation*}
Here for each $k$,  $S_{2k+1,2k}=\bigcap_{\varepsilon}(s_{\varepsilon})_* \bigl(H_*(K_{\varepsilon})\bigr)$, where $\varepsilon$ runs over sequences $\varepsilon$ of $\pm1$'s with $\varepsilon_k=1$ and $\varepsilon_{k+1}=0$, and $s_{\varepsilon}:K_{\varepsilon} \to K_{2k}\cap K_{2k+1}\subset P_{2k+1,2k}$ is the inclusion map. 
\end{Theorem C}

\begin{figure}
\begin{center}
\begin{tikzpicture}
\draw (0,0) -- (60: 1cm) (0,0) -- (320:1cm) (0,0) -- (300:1cm) 
 (130:1cm) -- (0,0) (150:1cm) -- (0,0) (230:1cm) -- (0,0);
\fill (0,0) circle (2pt);
\fill (60: 1cm) circle (2pt);
\fill (320:1cm) circle (2pt);
\fill (300:1cm) circle (2pt);
\fill (130:1cm) circle (2pt);
\fill (150:1cm) circle (2pt);
\fill (230:1cm) circle (2pt); 
\draw[dashed] (330:0.9cm) arc (330:360:0.8cm) arc (0:50:0.8cm);
\draw[dashed] (160:0.9cm) arc (160: 220:0.8cm);
\path (60: 1cm) node[above] {$K_{2k+2}$};
\path (320:1cm) node[right] {$K_{2m-2}$};
\path (300:1cm) node[below] {$K_{2m}$};
\path (130:1cm) node[above] {$K_1$};
\path (150:1cm) node[left] {$K_3$};
\path (230:1cm) node[left] {$K_{2k-1}$};
\path (50:0.9 cm) node[right] {$\gamma_{k+1}$};
\path (335:1.3cm) node[above] {$\gamma_{m-1}$};
\path (300:0.8cm) node[left] {$\gamma_m$};
\path (130:0.8cm) node[right] {$\gamma_1$};
\path (150:0.6cm) node[below] {$\gamma_2$};
\path (225:0.8cm) node[above] {$\gamma_k$};
\path (0,0) node[right] {$K_{\varepsilon}$};
\draw[->,>=stealth] (130:0.5cm) -- (130:0.49cm);
\draw[->,>=stealth] (150:0.5cm) -- (150:0.49cm);
\draw[->,>=stealth] (230:0.5cm) -- (230:0.49cm);
\draw[->,>=stealth] (60:0.5cm) -- (60:0.51cm);
\draw[->,>=stealth] (320:0.5cm) -- (320:0.51cm);
\draw[->,>=stealth] (300:0.5cm) -- (300:0.51cm);

\draw[->,>=stealth] (-1.7,0) -- (-2.7,0);
\path (-2.2,0) node[above] {$\iota_{\varepsilon}$};

\draw (-7.5,0.5) -- (-7.5,-0.5) (-7,0.5) -- (-7,-0.5) 
     (-6,0.5) -- (-6,-0.5)   (-5,0.5) -- (-5,-0.5) 
     (-4,0.5) -- (-4,-0.5) (-3,0.5) -- (-3,-0.5);
\fill  (-7.5,0.5) circle (2pt);
\fill (-7.5,-0.5) circle (2pt);
\fill (-7,0.5) circle (2pt);
\fill (-7,-0.5) circle (2pt);
\fill (-6,0.5) circle (2pt);
\fill (-6,-0.5) circle (2pt);
\fill (-5,0.5) circle (2pt);
\fill (-5,-0.5) circle (2pt);
\fill (-4,0.5) circle (2pt);
\fill (-4,-0.5) circle (2pt);
\fill (-3,0.5) circle (2pt);
\fill (-3,-0.5) circle (2pt);
\draw[->, >=stealth] (-7.5,-0.2) -- ++(0,-0.01);
\draw[->, >=stealth] (-7,-0.2) -- ++(0,-0.01);
\draw[->, >=stealth] (-6,-0.2) -- ++(0,-0.01);
\draw[->, >=stealth] (-5,-0.2) -- ++(0,-0.01);
\draw[->, >=stealth] (-4,-0.2) -- ++(0,-0.01);
\draw[->, >=stealth] (-3,-0.2) -- ++(0,-0.01);
\draw[dashed] (-4.2,0) -- (-4.8,0) (-6.2,0) -- (-6.8,0);
\path (-7.5,0.5) node[above] {$K_1$};
\path (-7.5,-0.5) node[below] {$K_2$};
\path (-7,0.5) node[above] {$K_3$};
\path (-7,-0.5) node[below] {$K_4$};
\path (-6,0.5) node[above] {$K_{2k-1}$};
\path (-6,-0.5) node[below] {$K_{2k}$};
\path (-5,0.5) node[above] {$K_{2k+1}$};
\path (-5,-0.5) node[below] {$K_{2k+2}$};
\path (-4,0.5) node[above] {$K_{2m-3}$};
\path (-4,-0.5) node[below] {$\ K_{2m-2}$};
\path (-3,0.5) node[above] {$\ \ K_{2m-1}$};
\path (-3,-0.5) node[below] {$K_{2m}$};
\path (-7.5,0.3) node[left] {$\gamma_1$};
\path (-7,0.3) node[left] {$\gamma_2$};
\path (-6,0.3) node[left] {$\gamma_k$};
\path (-5,0.3) node[left] {$\gamma_{k+1}$};
\path (-4,0.3) node[left] {$\gamma_{m-1}$};
\path (-3,0.3) node[left] {$\gamma_m$};

\end{tikzpicture}
\end{center}

\begin{center}
\begin{tikzpicture}
\draw[->,>=stealth] (1.7,0) -- (2.7,0);
\path (2.2,0) node[above] {$j_{\varepsilon}$};

\draw (3,0.5) -- (3,-0.5) (4,0.5) -- (4,-0.5) 
(6,0.5) -- (6,-0.5)  (7.5,0.5) -- (7.5,-0.5) 
(9,0.5) -- (9,-1.5);
\draw[dashed] (3.2,-0.1) -- (3.8,-0.1);
\draw[dashed] (6.2,-0.1) -- (7.3,-0.1);
\fill (3,0.5) circle (2pt);
\fill (3,-0.5) circle (2pt);
\fill (4,0.5) circle (2pt);
\fill (4,-0.5) circle (2pt);
\fill (6,0.5) circle (2pt);
\fill (6,-0.5) circle (2pt);
\fill (7.5,0.5) circle (2pt);
\fill (7.5,-0.5) circle (2pt);
\fill (9,0.5) circle (2pt);
\fill (9,-0.5) circle (2pt);
\fill (9,-1.5) circle (2pt);
\fill (5,0) circle (2pt);
\path (3,0.6) node[above] {$\eta_1$};
\path (4,0.6) node[above] {$\eta_{k-1}$};
\path (5,0.6) node[above] {$\eta_k$};
\path (6,0.6) node[above] {$\eta_{k+1}$};
\path (7.5,0.6) node[above] {$\eta_{m-1}$};
\path (9,0.6) node[above] {$\eta_m$};
\draw[->,>=stealth] (3,-0.2) -- ++(0,-0.01);
\draw[->,>=stealth] (4,-0.2) -- ++(0,-0.01);
\draw[->,>=stealth] (6,-0.2) -- ++(0,-0.01);
\draw[->,>=stealth] (7.5,-0.2) -- ++(0,-0.01);
\draw[->,>=stealth] (9,-0.2) -- ++(0,-0.01);
\draw[->,>=stealth] (9,-1.2) -- ++(0,-0.01);
\path (3,0.5) node[right] {$K_3$};
\path (3,-0.5) node[right] {$K_{\varepsilon}$};
\path (4,0.5) node[right] {$K_{2k-1}$};
\path (4,-0.5) node[right] {$K_{\varepsilon}$};
\path (5,0) node[right] {$K_{\varepsilon}$};
\path (6,0.5) node[right] {$K_{\varepsilon}$};
\path (6,-0.5) node[right] {$K_{2k+2}$};
\path (7.5,0.5) node[right] {$K_{\varepsilon}$};
\path (7.5,-0.5) node[right] {$K_{2m-2}$};
\path (9,0.5) node[right] {$K_1$};
\path (9,-0.5) node[right] {$K_{\varepsilon}$};
\path (9,-1.5) node[right] {$K_{2m}$};
\path (3,0.1) node[right] {$\gamma_2$};
\path (4,0.1) node[right] {$\gamma_k$};
\path (6,0.1) node[right] {$\gamma_{k+1}$};
\path (7.5,0.1) node[right] {$\gamma_{m-1}$};
\path (9,0.1) node[right] {$\gamma_1$};
\path (9,-1) node[right] {$\gamma_m$};

\draw node[text width=11cm]  at (5,-4)
{\textsc{Figure 4.} A simultaneous saddle interaction of $m$ incoming open strings $\gamma_1$,$\dots$,$\gamma_m$ when the interaction parameter $\vec{t}\in [0,1]^m$ is a corner value $\varepsilon=(1,\dots,1,0,\dots,0)$ with the last $1$ is at the $k$-th position. Here the first $k$ incoming open strings meet at their head vertices and the remaining $m-k$ open strings meet at their tail vertices. The result of the interaction is $m$ outgoing open strings $\eta_1,\dots,\eta_m$. Note that the $k$-th outgoing open string $\eta_k$ is a constant open string contained in the submanifold $K_{\varepsilon}$. This observation leads to the proof of Theorem C (Theorem 5.1). See \S5 for further details. };
\end{tikzpicture}
\end{center}
\end{figure}

For example, suppose $\Sigma$ has two incoming open strings and two outgoing open strings, arranged alternately. Let $K_1,K_2,K_3,K_4$ be oriented closed submanifolds of $M$ labeling four free boundary components of $\Sigma$ so that $P_{K_1K_2}\times P_{K_3K_4}$ denotes the space of incoming open strings, and $P_{K_3K_2}\times P_{K_1K_4}$ denotes the space of outgoing open strings. The open-closed cobordism $\Sigma$ describes the simplest saddle interaction. See Figure 5. The associated saddle string operation factors through the homology group of constant open strings as follows (Corollary \ref{m=2 example}).  
\begin{equation*}
\xymatrix{
H_*(P_{K_1K_2})\otimes H_*(P_{K_3K_4}) 
\ar[r]^{\mu_{\Sigma}\ \ \ \ } \ar[dr] 
& H_*(P_{K_3K_2})\otimes H_*(P_{K_1K_4}) \\
& H_*(K_3\cap K_2)\otimes H_*(K_1\cap K_4) \ar[u]_{s_*\otimes s_*}
}
\end{equation*}
If $K_2\cap K_3=\emptyset$ or $K_1\cap K_4=\emptyset$, then we see that the above saddle string operation vanishes. A general vanishing property of this type is given in Corollary \ref{vanishing of string operation}. For the case of the saddle interaction of three incoming open strings, see Corollary \ref{m=3 example}. Theorem C is proved in Theorem 5.1 in \S 5.1. 

The deformation process in the proof of Theorem C through simultaneous saddle interactions can also be used to prove that the string operation $\mu_{\Sigma}$ is independent of its half-pair-of-pants decomposition.

\begin{figure}
\begin{center}
\begin{tikzpicture}[scale=0.3, auto] 
\path (-0.8,3) coordinate (1) node[left] {$K_1$};
\path (-3.42,4) coordinate (2) node[above] {$K_2$};
\path (-0.58,-4) coordinate (3) node[below] {$K_4$};
\path (-3,-3) coordinate (4) node[below] {$K_3$};
\path (4,3) coordinate (5) node[above] {$K_1$};
\path (1.58,4) coordinate (6) node[above] {$K_2$};
\path (4.42,-4) coordinate (7) node[below] {$K_4$};
\path (2,-3) coordinate (8) node[right] {$K_3$};
\path (-0.79,-3.25) coordinate (9);
\path (1.8,3.25) coordinate (10);
\path (2.53,0) coordinate (11);
\path (3.48,0) coordinate (12);
\draw (-2,-1) arc (270:360: 2 cm and 1 cm) 
arc (0:90:2 cm and 1 cm);
\draw[ultra thick] (1) parabola bend (-2,1) (2);
\draw[->, >=stealth, ultra thick] (-2.3,1) -- ++(-0.1,0);
\draw (2) to [out=350, in=200]  (6) ;
\draw (1) to [out=10, in=170]  (5);
\draw[ultra thick] (4) parabola bend (-2,-1) (3);
\draw[->, >=stealth, ultra thick] (-1.7,-1) -- ++(0.1,0);
\draw (3) to [out=10,in=170] (7);
\draw[ultra thick] (7) to [out=110,in=250]  (5);
\path (2,-3) coordinate (8);
\draw (4) to [out=350,in=175] (9);
\draw[densely dashed] (9) to [out=355,in=190] (8);
\draw[densely dashed, ultra thick] (8) to [out=70,in=290] (10);
\draw[ultra thick] (10) -- (6);
\draw (0,0) arc (180:253:5 cm and 0.7 cm);
\draw[densely dashed] (0,0) arc (180:120:5 cm and 0.7 cm);
\draw[->, >=stealth, ultra thick] (11) -- (2.53,0.2);
\draw[->, >=stealth, ultra thick] (3.48,0.01) -- (12);

\fill (1) circle (4pt);
\fill (2) circle (4pt);
\fill (3) circle (4pt);
\fill (4) circle (4pt);
\fill (5) circle (4pt);
\fill (6) circle (4pt);
\fill (7) circle (4pt);
\fill (8) circle (4pt);

\draw node[text width=11cm] at (0,-8) 
{\textsc{Figure 5}. The saddle interaction of open strings has the following form and has values in constant homology classes (see Corollary 5.3):
$H_*(P_{K_1K_2})\otimes H_*(P_{K_3K_4}) \rightarrow H_*(P_{K_3K_2})\otimes H_*(P_{K_1K_4})$.};
\end{tikzpicture}
\end{center}
\end{figure}

\begin{Theorem D} Let $\Sigma$ be an open-closed cobordism homeomorphic to a disc. Then the associated string operation $\mu_{\Sigma}$ is independent of the half-pair-of-pants decomposition of $\Sigma$. Furthermore, the same string operation can be computed by a simultaneous saddle interaction of incoming open strings at the same point \textup{(}see Figure 3 and Figure 4\textup{)}, which may be internal points or end points of open strings. 
\end{Theorem D}

The organization of the paper is as follows. In \S 2, we introduce sewing diagrams which describe instructions on how open string products and coproducts are to be composed to produce string operations. In \S 3, saddle interaction diagrams are introduced to describe more general interactions of open strings at their internal points and discuss their degenerated simultaneous saddle interactions. This leads to the proof of Theorem D given as Theorem \ref{independence of decomposition}. In \S 4, we show that iterated coproducts have values in homology classes of constant open strings. We also discuss its relation to open string products with homology classes of constant open strings. These lead to proofs of Theorem A and Theorem B given in Proposition \ref{only outgoing open strings} and in Theorem \ref{image of open coproducts}. Finally in \S 5, we describe a qualitative behavior of string operations associated to disc cobordisms by considering various degenerated simultaneous saddle interactions and prove Theorem C.

\bigskip

\section{Sewing diagrams and open string interactions}

Let $\Sigma$ be an open-closed cobordism homeomorphic to a disc. For each choice of its decomposition into half-pair-of-pants, we can compute the associated string operation by composing open string products and coproducts according to the decomposition. The decomposition data can be conveniently encoded into sewing diagrams and their decompositions. In \S 3 we will discuss relations among different sewing diagrams and show that the string operation $\mu_{\Sigma}$ associated to a disc cobordism is independent of the decomposition, up to sign. The main ingredient of our proof is open string interactions at internal points, which we call saddle interactions. 

\subsection{Sewing diagrams} 

We first define various objects related to sewing diagrams. 

A \emph{sewing graph} $H$ is a finite oriented forest whose vertices have degree at most 3, and whose vertices of degree $2$ or $3$ are neither sinks nor sources. 

Let $G$ be a disjoint union of oriented edges and let $H$ be a sewing graph. The graph $G$ can be regarded as a trivial sewing graph. A \emph{sewing morphism} $f: G \rightarrow H$ from $G$ to $H$ is an orientation preserving surjective map which can be decomposed into a sequence of maps 
\begin{equation}\label{sewing morphism decomposition}
G \xrightarrow{f_1} H_1  \xrightarrow{f_2}  H_2  \xrightarrow{f_3}  {\dotsb}
\xrightarrow{f_r} H_r=H,
\end{equation}
where $H_0=G, H_1, ,H_2,\dotsc,H_r$ are sewing graphs and $f_1,f_2,\dotsc,f_r$ are maps of the following types:
\begin{enumerate}
\item[(I)]  $f_i$ identifies two end vertices from different connected components, where one of them is a head vertex and the other is a tail vertex. 
\item[(II)]  $f_i$ inserts a vertex into the interior of an edge.
\item[(III)]  $f_i$ identifies an end vertex and a degree $2$ vertex arising as in (II) in different components. 
\end{enumerate} 

Thus, for a sewing morphism $f:G \rightarrow H$, the inverse image of a degree $2$ vertex is either two end vertices from distinct components, or an internal point of an edge. In the first case, the degree $2$ vertex is called of type (I), and in the latter case, of type (II). The inverse image in $G$ of a degree $3$ vertex in $H$ always consists of an internal point of an edge and an end vertex of another edge.

A sewing morphism $f:G \rightarrow H$ can be canonically extended to a \emph{sewing diagram} of the form 
\begin{equation*}
G \xrightarrow{f} H \xleftarrow{f'} G',
\end{equation*}
where $f': G' \rightarrow H$ is another sewing morphism, by performing the following operations to vertices $v$ of $H$ of degree $2$ and $3$. When the vertex $v$ has degree $3$, let $e_1,e_2,e_3$ be edges in $H$ meeting at $v$ such that $e_2$ and $e_3$ come from the same edge in $G$, and $e_2$ is directed toward $v$ and $e_3$ is directed away from $v$.
\begin{enumerate}
\item If $v$ is a degree $2$ vertex of type (I), then erase the vertex from $H$ to form a new oriented edge. The orientation of the new edge is well defined since the orientation of $H$ is concordant at $v$.  

\begin{center}
\begin{tikzpicture}
\draw[thick] (0,0) -- ++(0.5,0);
\draw[thick] (0.8,0) -- ++(0.5,0);
\fill (0.5,0) circle (2pt);
\fill (0.8,0) circle (2pt);
\draw[->, >=stealth, thick] (0.25,0) -- ++(0.1,0);
\draw[->,>=stealth, thick] (1.05,0) -- ++(0.1,0);
\draw[->,>=stealth] (1.65,0) -- ++(0.5,0);
\draw[thick] (2.5,0) -- ++(1,0);
\path (3,0) coordinate node[above] {$v$};
\fill (3,0) circle (2pt);
\draw[->,>=stealth, thick] (2.75,0) -- ++(0.1,0);
\draw[->,>=stealth, thick] (3.25,0) -- ++(0.1,0);
\draw[->,>=stealth] (4.2,0) -- ++(-0.5,0);
\draw[thick] (4.55,0) -- ++(1,0);
\draw[->,>=stealth, thick] (5,0) -- ++(0.1,0);
\end{tikzpicture}
\end{center}

\item If $v$ is a degree $2$ vertex of type (II), then split $v$ into two end vertices carrying the same label as $v$, where edges adjacent to new vertices have the same orientation as the original edge in $G$. 

\begin{center} 
\begin{tikzpicture}
\draw[thick] (0,0) -- (1,0);
\draw[->,>=stealth, thick] (0.5,0) -- ++(0.1,0);
\draw[->,>=stealth] (1.3,0) -- ++(0.5,0);
\draw[thick] (2.1,0) -- ++(1,0);
\fill (2.6,0) circle (2pt);
\path (2.6,0) node[above] {$v$};
\draw[->,>=stealth, thick] (2.3,0) -- ++(0.1,0);
\draw[->,>=stealth, thick]  (2.9,0) -- ++(0.1,0);
\draw[->,>=stealth] (3.9,0) -- ++(-0.5,0);
\draw[thick] (4.2,0) -- (4.7,0);
\fill (4.7,0) circle (2pt);
\draw[thick] (5,0) -- (5.5,0);
\fill (5,0) circle (2pt);
\draw[->,>=stealth, thick] (4.4,0) -- ++(0.1,0);
\draw[->,>=stealth, thick]  (5.2,0) -- ++(0.1,0);
\end{tikzpicture}
\end{center}

\item If the vertex $v$ is of degree $3$ and the edge $e_1$ is directed toward $v$, then cut $H$ at $v$ so that the edge $e_2$ ends at $v$, and  edges $e_1$ and $e_3$ are joined in this order and form a new oriented edge by erasing the vertex $v$. 

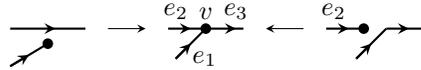
\begin{center} 
\begin{tikzpicture}
\draw[thick] (0,0) -- (1,0);
\draw[thick] (0,-0.5) -- (0.5,-0.2);
\fill (0.5,-0.2) circle (2pt);
\draw[->, >=stealth, thick] (0,-0.5) -- (0.25,-0.35);
\draw[->,>=stealth, thick] (0.5,0) -- ++(0.1,0);
\draw[->,>=stealth] (1.3,0) -- ++(0.5,0);
\draw[thick] (2.1,0) -- ++(1,0);
\path (2.2,0) node[above] {$e_2$};
\path (3,0) node[above] {$e_3$};
\fill (2.6,0) circle (2pt);
\path (2.6,0) node[above] {$v$};
\draw[thick] (2.2,-0.4) -- (2.6,0);
\path (2.3,-0.4) node[right] {$e_1$};
\draw[->, >=stealth, thick] (2.2,-0.4) -- (2.4,-0.2);
\draw[->,>=stealth, thick] (2.3,0) -- ++(0.1,0);
\draw[->,>=stealth, thick]  (2.9,0) -- ++(0.1,0);
\draw[->,>=stealth] (3.9,0) -- ++(-0.5,0);
\draw[thick] (4.2,0) -- (4.7,0);
\fill (4.7,0) circle (2pt);
\path (4.3,0) node[above] {$e_2$};
\draw[thick] (5,0) -- (5.5,0);
\draw[thick] (4.6,-0.4) -- (5,0);
\draw[->, >=stealth, thick] (4.6,-0.4) -- (4.8, -0.2);
\draw[->,>=stealth, thick] (4.4,0) -- ++(0.1,0);
\draw[->,>=stealth, thick]  (5.2,0) -- ++(0.1,0);
\end{tikzpicture}
\end{center}

\item If the vertex $v$ is of degree $3$ and the edge $e_1$ is directed away from $v$, then cut $H$ at $v$ so that $e_3$ starts with the vertex $v$, and edges $e_2$ and $e_1$ are joined in this order to form a new oriented edge by erasing the vertex $v$. 

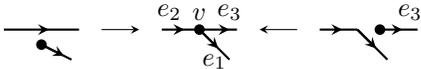
\begin{center} 
\begin{tikzpicture}
\draw[thick] (0,0) -- (1,0);
\draw[->,>=stealth, thick] (0.5,0) -- ++(0.1,0);
\draw[thick] (0.5,-0.2) -- (0.9,-0.4);
\fill (0.5,-0.2) circle (2pt);
\draw[->, >=stealth, thick] (0.5,-0.2) -- (0.8, -0.35);
\draw[->,>=stealth] (1.3,0) -- ++(0.5,0);
\draw[thick] (2.1,0) -- ++(1,0);
\path (2.2,0) node[above] {$e_2$};
\path (3,0) node[above] {$e_3$};
\fill (2.6,0) circle (2pt);
\path (2.6,0) node[above] {$v$};
\draw[thick] (2.6,0) -- (3,-0.4);
\path (2.8,-0.2) node[below] {$e_1$};
\draw[->, >=stealth, thick] (2.6,0) -- (2.85, -0.25);
\draw[->,>=stealth, thick] (2.3,0) -- ++(0.1,0);
\draw[->,>=stealth, thick]  (2.9,0) -- ++(0.1,0);
\draw[->,>=stealth] (3.9,0) -- ++(-0.5,0);
\draw[thick] (4.2,0) -- (4.7,0);
\draw[thick] (4.7, 0) -- ++(0.4,-0.4);
\draw[->, >=stealth, thick] (4.7,0) -- ++(0.25, -0.25);
\draw[thick] (5,0) -- (5.5,0);
\fill (5,0) circle (2pt);
\path (5.4,0) node[above] {$e_3$};
\draw[->,>=stealth, thick] (4.4,0) -- ++(0.1,0);
\draw[->,>=stealth, thick]  (5.2,0) -- ++(0.1,0);
\end{tikzpicture}
\end{center}
\end{enumerate} 
Let $G'$ be the graph resulting from these operations performed on $H$. Since $G'$ is a directed graph with only degree $1$ vertices, it must be a disjoint union of oriented edges, and the map $f':G'\rightarrow H$, undoing the above operations on $H$, is a sewing morphism. Note that $(f')'=f$. 

For a disjoint union $G$ of oriented edges, a trivial sewing diagram is a diagram of the form $G \xrightarrow{f} G \xleftarrow{f'} G$ where both $f$ and $f'$ are identity maps. 

When we have two sewing diagrams $G_i \xrightarrow{f_i} H_i \xleftarrow{f_i'} G'_i$ for $i=1,2$, their union is a sewing diagram of the form 
\begin{equation*}
G_1+G_2 \xrightarrow{f_1+f_2} H_1+H_2 \xleftarrow{f_1'+f_2'} G_1'+G_2', 
\end{equation*}
where $+$ denotes a disjoint union of two graphs. 

A sewing diagram $G \xrightarrow{f} H \xleftarrow{f'} G'$ is called type (I) if $f$ is a type (I) morphism identifying exactly two end vertices from different edges. It is called type (II) if $f$ is a type (II) morphism inserting a vertex to the interior of an edge. Type (I) sewing diagrams correspond to open string products, and type (II) sewing diagrams correspond to open string coproducts. See \S 2.3. 

Let $G\xrightarrow{f} H \xleftarrow{f'} G'$ be a sewing diagram. When vertices of $H$ carry labels, we can consistently label vertices of $G$ in such a way that a vertex $v$ in $G$ carry the same label 
as the image vertex $f(v)$ in $H$. Similarly, we can label vertices in $G'$. When vertices of a sewing diagram are labeled in this way, we call it a labeled sewing diagram. Thus, in a type (I) labeled sewing diagram, two vertices identified by $f:G \to H$ carry the same label. Similarly for a type (II) labeled sewing diagram. 

Let $G \xrightarrow{f} H \xleftarrow{f'} G'$ be a sewing diagram, where both $G$ and $G'$ are disjoint union of oriented edges. The only topological invariant of $G$ and $G'$ is their number of edges which we denote by $e$ and $e'$. Let $p_1$ and $p_{2}$ denote the number of vertices of degree $2$ of type (I) and (II), respectively, in $H$, and let $p_{3}$ denote the number of vertices of degree $3$ in $H$.  Then the number of edges $e(H)$ of $H$ is given by $e(H)=e+p_2+p_{3}$ on the one hand since $p_2+p_3$ internal points of $G$ become vertices, and $e(H)=e'+p_1+p_{3}$ on the other hand. Hence the number of edges of $G$ and $G'$ are related by $e'=e+p_1-p_2$. In particular, if $H$ has no degree $2$ vertices, then $G$ and $G'$ have the same number of edges.  

Here is a simplest example of a sewing diagram. 

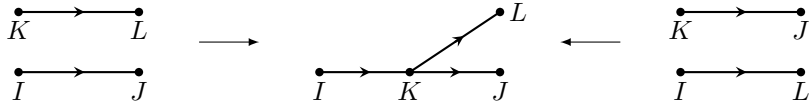
\begin{center}
\begin{tikzpicture}[scale=0.8]
\path (0,0) coordinate (I1) node[below] {$I$};
\fill (I1) circle (2pt) ; 
\path (2,0) coordinate (J1) node[below] {$J$};
\fill (J1) circle (2pt) ; 
\path (0,1) coordinate (K1) node[below] {$K$};
\fill (K1) circle (2pt) ; 
\path (2,1) coordinate (L1) node[below] {$L$};
\fill (L1) circle (2pt) ; 
\path (5,0) coordinate (I2) node[below] {$I$};
\fill (I2) circle (2pt) ; 
\path (6.5,0) coordinate (K2) node[below] {$K$};
\fill (K2) circle (2pt) ; 
\path (8,0) coordinate (J2) node[below] {$J$};
\fill (J2) circle (2pt) ; 
\path (8,1) coordinate (L2) node[right] {$L$};
\fill (L2) circle (2pt) ; 
\path (11,1) coordinate (K3) node[below] {$K$};
\fill (K3) circle (2pt) ; 
\path (13,1) coordinate (J3) node[below] {$J$};
\fill (J3) circle (2pt) ; 
\path (11,0) coordinate (I3) node[below] {$I$};
\fill (I3) circle (2pt) ; 
\path (13,0) coordinate (L3) node[below] {$L$};
\fill (L3) circle (2pt) ; 

\path (3,0.5) coordinate (a);
\path (4,0.5) coordinate (b) ;
\path (9,0.5) coordinate (c);
\path (10,0.5) coordinate (d);

\draw[thick] (I1)-- (J1) (K1) -- (L1) (I2) -- (J2) (K2) -- (L2) (I3) -- (L3) (K3) -- (J3);
\draw[->, >=stealth, thick] (1,0) -- ++(0.1,0);
\draw[->, >=stealth, thick] (1,1) -- ++(0.1,0);
\draw[->, >=stealth, thick] (5.75,0) -- ++(0.1,0);
\draw[->, >=stealth, thick] (7.25,0) -- ++(0.1,0);
\draw[->, >=stealth, thick] (7.25,0.5) -- ++(0.15,0.1);
\draw[->, >=stealth, thick] (12,0) -- ++(0.1,0);
\draw[->, >=stealth, thick] (12,1) -- ++(0.1,0);

\draw[-latex] (a) -- (b);
\draw[-latex] (d) -- (c);

\end{tikzpicture}
\end{center}

This sewing diagram represents a composition of an open string coproduct of $P_{IJ}$ along $K$ followed by an open string product of $P_{IK}$ and $P_{KL}$ along $K$:
\begin{multline*}
H_*(P_{IJ})\otimes H_*(P_{KL}) \xrightarrow{\varphi_K\otimes1} 
H_*(P_{IK})\otimes H_*(P_{KJ})\otimes H_*(P_{KL}) \\
\xrightarrow{1\otimes T} H_*(P_{IK})\otimes H_*(P_{KL})\otimes H_*(P_{KJ}) \xrightarrow{\mu_K\otimes 1} H_*(P_{IL})\otimes H_*(P_{KJ})
\end{multline*}
See Figure 6 for the corresponding open cobordism. 
Any degree $3$ vertex in a sewing diagram represents a composition of this type, namely a saddle operation.

\begin{figure}
\begin{center}
\begin{tikzpicture}[scale=0.7]
\draw (0,0) -- (0,1) to[out=0, in=225] 
(2,2) to[out=45, in=180] (4,3) -- (6,3)
-- (6,2) -- (4,2) to[out=180, in=45] (3,1.5) 
[rounded corners] -- (2.5,1) [rounded corners] -- (3.5,1) 
--(3,1.5);
\draw (0,0) -- (6,0) -- (6,1) to[out=180,in=315] (4,2);
\draw (2,2) -- (0,2) -- (0,3) -- (2,3) to[out=360,in=150] 
(3,2.8);
\draw[thick] (0,0) -- (0,1) (0,2) -- (0,3) (6,0) -- (6,1) (6,2) -- (6,3);
\draw[->,>=stealth, thick] (0,0.6) -- ++(0,0.1);
\draw[->,>=stealth, thick] (6,0.6) -- ++(0,0.1);
\draw[->,>=stealth, thick] (0,2.6) -- ++(0,0.1);
\draw[->,>=stealth, thick] (6,2.6) -- ++(0,0.1);
\path (0,0) node[left] {$I$};
\path (0,1) node[left] {$J$};
\path (0,2) node[left] {$K$};
\path (0,3) node[left] {$L$};
\path (6,0) node[right] {$I$};
\path (6,1) node[right] {$L$};
\path (6,2) node[right] {$K$};
\path (6,3) node[right] {$J$};

\path node[text width=11cm] at (3,-2)
{\textsc{Figure 6.} An open-closed cobordism representing an open string coproduct followed by an open string product. This is the same as the open-closed cobordism for the saddle operation in Figure 5. };
\end{tikzpicture}
\end{center}
\end{figure}

\subsection{Composition and decomposition of sewing diagrams}

Next we will show that any labeled sewing diagram is a composition of type (I) labeled sewing diagrams and type (II) labeled sewing diagrams in some order. 

First we describe how to compose an arbitrary labeled sewing diagram with a type (I) labeled sewing diagram. Let 
$G \xrightarrow{f} H \xleftarrow{f'} G'$ be a labeled sewing diagram. Suppose $G'$ has two vertices $v_1$ and $v_2$ from different edges carrying the same label such that $f'(v_1)$ and $f'(v_2)$ are distinct vertices in $H$, and identifying $f'(v_1)$ and $f'(v_2)$ in $H$ gives a map $g''$ of type (I) or (III) with respect to $f$. Thus either both $f'(v_1)$ and $f'(v_2)$ are end vertices from different components of $H$ such that exactly one of them is a head vertex, or one of them is an end vertex and the other is a vertex of type (II) with respect to $f$ from another component of $H$. Let $G' \xrightarrow{g} H' \xleftarrow{g'} G''$ be a type (I) sewing diagram in which $g$ identifies vertices $v_1$ and $v_2$ into a vertex $v'$. Since $v_1$ and $v_2$ carry the same label, $H'$ and $G''$ can be consistently labeled. Now we have the following diagram $(a)$ in \eqref{composition of sewing diagrams} in which $g''$ is an identification map of two vertices $f'(v_1)$ and $f'(v_2)$ into a vertex $v''$, and $f''$ is induced from $f'$ by mapping $v'$ in $H'$ to $v''$ in $H''$.
\begin{equation}\label{composition of sewing diagrams} 
(a) \quad
\xymatrix{
G \ar[dr]_{g''\circ f} \ar[r]^f & H \ar[d]^{g''} 
& G' \ar[l]_{f'}  \ar[d]^g \\
& H'' & H'\ar[l]_{f''} \\
& & G'' \ar[ul]^{f''\circ g'} \ar[u]_{g'} }
\qquad\quad
(b) \quad
\xymatrix{
G \ar[dr]_{h''\circ f} \ar[r]^f & H \ar[d]^{h''} 
& G' \ar[l]_{f'}  \ar[d]^h \\
& H'' & H'\ar[l]_{f''} \\
& & G'' \ar[ul]^{f''\circ h'} \ar[u]_{h'} }
\end{equation}
Since $H''$ has vertices of degree at most $3$ and $g'': H \rightarrow H''$ is a type (I) map (if both vertices $f'(v_1)$ and $f'(v_2)$ have degree 1) or type (III) map (if exactly one of the vertices $f'(v_1)$ or $f'(v_2)$ has degree 2) with respect to $f$, the composition $g''\circ f: G \rightarrow H''$ is a sewing morphism, and the other sewing morphism $f''\circ g': G'' \rightarrow H''$ is induced from it. The resulting labeled sewing diagram $T'': G \rightarrow H'' \leftarrow G''$ is defined to be the composition of the labeled sewing diagrams $T: G \xrightarrow{f} H \xleftarrow{f'} G'$ and $T': G' \xrightarrow{g} H' \xleftarrow{g'} G''$, and we write $T''=T'\circ T$. 

Next, we describe a composition of an arbitrary labeled sewing diagram with a type (II) labeled sewing diagram. For a sewing diagram $G \xrightarrow{f} H \xleftarrow{f'} G'$, let $p$ be an internal point of an edge $e$ in $G'$ such that $f'(p)$ is not a vertex in $H$. In this case, let $G' \xrightarrow{h} H' \xleftarrow{h'} G''$ be a labeled sewing diagram of type (II) inserting a new vertex $v'$ at the internal point $p$ of an edge $e$ in $G'$. Let $h'': H \rightarrow H''$ be a map inserting a vertex $v''$ at the internal point $f'(p)$ in $H$. We have the diagram $(b)$ in \eqref{composition of sewing diagrams} in which $f''$ maps the vertex $v'$ to the vertex $v''$. 
Again, since $H''$ has vertices of degree at most $3$, and $h''$ is a map of type (II), the composition $h''\circ f: G \rightarrow H''$ is a sewing morphism. The resulting labeled sewing diagram $T'': G \rightarrow H'' \leftarrow G''$ is by definition the composition of the labeled sewing diagrams $T: G \rightarrow H \leftarrow G'$ and $T': G' \rightarrow H' \leftarrow G''$. We write $T''=T'\circ T$. 

Let $T: G \xrightarrow{f} H \xleftarrow{f'} G'$ be a labeled sewing diagram. For each choice of a decomposition of the sewing morphism $f$ into maps of types (I), (III), or type (II) as in \eqref{sewing morphism decomposition}, there corresponds a sequence of type (I) or type (II) sewing diagrams $T_1,T_2,\dotsc, T_r$ such that $T$ is a composition of these elementary sewing diagrams: 
\begin{equation}\label{sewing diagram decomposition}
T=T_r\circ \dotsb \circ T_2\circ T_1.
\end{equation}
Conversely, any such decomposition of a given sewing diagram into type (I) and type (II) sewing diagrams gives a decomposition of $f: G \rightarrow H$ in the form \eqref{sewing morphism decomposition}. In fact, these two types of decompositions are in bijective correspondence. 

To each decomposition \eqref{sewing diagram decomposition} of a labeled sewing diagram, we can associate an open cobordism from open strings of type $G$ to open strings of type $G'$ with labeled free boundaries as follows. 

Given a type (I) sewing diagram $G \xrightarrow{\ g} H \xleftarrow{g'} G'$, let $e_1$ and $e_2$ be two distinct oriented edges in $G$ which are joined by $g$ into an oriented edge $e_3$ in $G'$. To this type (I) sewing diagram, we associate an open cobordism consisting of an open string product cobordism with two incoming open strings $e_1, e_2$ and one outgoing open string $e_3$, and trivial open cobordisms for all the other edges in $G$. Similarly, for a type (II) sewing diagram, we associate an open cobordism consisting of an appropriate open string coproduct cobordism and trivial open cobordisms. The labeling of free boundaries is done in an obvious way. 

Given a decomposition of a sewing diagram $T: G\xrightarrow{f} H \xleftarrow{f'} G'$ in the form \eqref{sewing diagram decomposition}, we sew open string cobordisms corresponding to $T_1, T_2, \dotsc, T_r$ in this order. Since type (I) sewing diagrams identify end vertices from different connected components of sewing graphs and type (II) sewing diagrams split edges, the resulting open string cobordism is homeomorphic to a disjoint union of discs, where the number of connected components of the open cobordism is the same as the number of connected components of the sewing graph $H$. 

Thus, when $H$ is connected, every choice of the decomposition of the sewing morphism $f:G \rightarrow H$ gives rise to an open cobordism homeomorphic to a disc equipped with a decomposition into half-pair-of-pants which are ordered according to the sequence of the decomposition of $f:G \rightarrow H$. 

Every open cobordism homeomorphic to a disc can be canonically decomposed into three parts. We first combine adjacent incoming open strings if any, then the resulting incoming open strings interact and they are recombined into outgoing open strings of equal number, then finally some of the outgoing open strings are split into several outgoing open strings.
This decomposition corresponds to a \emph{canonical decomposition} of a given sewing diagram $T: G\xrightarrow{f} H \xleftarrow{f'} G'$ in the form $T=T_{2}\circ T_3\circ T_1$ constructed as follows. Let $H_1$ be a sewing graph obtained from $G$ by identifying vertices in $G$ corresponding to type (I) degree 2 vertices in $H$. We have a sewing diagram $T_1: G \rightarrow H_1 \leftarrow G''$. Let $H_3$ be a sewing graph with only vertices of degree 1 and 3 obtained from $H$ by erasing all degree $2$ vertices. We have a sewing diagram $T_3: G'' \rightarrow H_3 \leftarrow G'''$. Finally, let $H_2$ be a sewing graph obtained from $G'$ by identifying vertices in $G'$ coming from type (II) degree $2$ vertices in $H$. We have a sewing diagram $T_2: G''' \rightarrow H_2 \leftarrow G'$. Each connected components of open string cobordism corresponding to $T_1$ is homeomorphic to a disc and has exactly one outgoing open string, and similarly, each connected component of open string cobordism for $T_2$ is homeomorphic to a disc and has exactly one incoming open string. Each connected component of the open string cobordism associated to $T_3$ is homeomorphic to a disc, and incoming open strings and outgoing open strings alternate along its boundary. Consequently, there are the same number of incoming and outgoing open strings on each such component. We state this canonical decomposition formally. 

\begin{prop} \label{canonical decomposition of sewing diagrams} Every labeled sewing diagram $T$ admits a canonical decomposition of the form $T=T_2\circ T_3\circ T_1$, where the labeled sewing graph $H_i$ associated to $T_i$ for $i=1,2,3$ has the following properties. 
\begin{enumerate}
\item $H_1$ is a disjoint union of path graphs with only type \textup{(I)} degree $2$ vertices. 
\item $H_2$ is a disjoint union of path graphs with only type \textup{(II)} degree $2$ vertices. 
\item $H_3$ is a forest with only vertices of degree $1$ or $3$. 
\end{enumerate} 
\end{prop}

\subsection{String operations associated to sewing diagrams} 

We can associate homology operations to sewing diagrams labeled by a family of oriented closed submanifolds of $M$. These operations are TQFT string operations of the associated open string cobordisms. 

Let $\{K_{\lambda}\}_{\lambda\in\Lambda}$ be a family of closed oriented submanifolds of $M$, so called D-branes. For a labeled graph $H$ with vertex labeling given by $\ell: V(H) \rightarrow \Lambda$, let 
\begin{equation*}
\textup{Map}_*(H,M)=\{\gamma: H \rightarrow M \mid \gamma(v)\in K_{\ell(v)} \textup{ for all } v\in V(H)\}
\end{equation*}
be the restricted mapping space consisting of continuous maps from $H$ (regarded as a topological space) to $M$ such that vertices of $H$ are mapped into submanifolds specified by their labels. Thus, for example, 
\begin{equation*}
\textup{Map}_*(V(H), M)\cong \prod_{v\in V(H)} K_{\ell(v)}
\end{equation*}
is an orientable manifold, and only when we specify an ordering of vertices of $H$, it becomes an oriented manifold. 

A labeled sewing diagram $T: G \xrightarrow{f} H \xleftarrow{f'} G'$  labeled by $\{K_{\lambda}\}_{\lambda\in\Lambda}$ induces the following pull-back diagram of fibrations: 
\begin{equation*} 
\begin{CD}
\textup{Map}_*(G, M) @<{f}<< \textup{Map}_*(H,M) @>{f'}>> 
\textup{Map}_*(G',M) \\
@V{p}VV   @V{p}VV   @V{p}VV   \\
\textup{Map}_*\bigl(f^{-1}\bigl(V(H)\bigr),M\bigr) 
@<{\bar{f}}<<  \textup{Map}_*\bigl(V(H),M\bigr) 
@>{\bar{f'}}>>  
\textup{Map}_*\bigl({f'}^{-1}\bigl(V(H)\bigr),M\bigr)
\end{CD}
\end{equation*} 
Here, we use the same notation $f$ and $f'$ for induced maps on mapping spaces. Note that $f^{-1}\bigl(V(H)\bigr)$ can contain internal points of $G$, as well as vertices of $G$. Similarly for ${f'}^{-1}\bigl(V(H)\bigr)$. 

Let $K$ and $K'$ be oriented closed submanifolds in the above family of labels, and let $P_{KK'}$ be the open string space of continuous paths from points in $K$ to points in $K'$. Then, since $G$ is a disjoint union of edges with vertex labels, the mapping space $\textup{Map}_*(G,M)$ is homeomorphic to a product of open string spaces $P_{KK'}$ for some $K,K'$. 

Since smooth manifolds in the bottom line are orientable finite dimensional manifolds, the smooth normal bundle to the embedding $\bar{f}$ is also orientable. Consequently, we have a Thom class of the normal bundle, well defined up to a sign. This sign ambiguity is due to the lack of preferred ordering of elements in the sets $V(H)$ and $f^{-1}(V(H))$. The Thom class  allows us to define a transfer map $f_!$ determined up to sign, and with a degree shift: 
\begin{equation*}
f_!: H_*\bigl(\textup{Map}_*(G,M)\bigr) \rightarrow 
H_*\bigl(\textup{Map}_*(H,M)\bigr) .
\end{equation*}
Let $V_1$ and $V_2$ be the set of vertices of $H$ of degree $2$ of type (I) and type (II), respectively, and let $V_3$ be the set of degree $3$ vertices of $H$.  The map $\bar{f}$ between base manifolds is a product of the following maps; 
\begin{align*} 
K_{\ell(v_1)}\times K_{\ell(v_1)} &\longleftarrow K_{\ell(v_1)} \textup{ \ for \ } v_1\in V_1, \\
M &\longleftarrow  K_{\ell(v_2)} \textup{ \ for \ } v_2\in V_2, \\
M\times K_{\ell(v_3)} &\longleftarrow  K_{\ell(v_3)} 
\textup{ \ for \ } v_3\in V_3.
\end{align*}
The homology operation associated to a labeled sewing diagram 
$T: G \xrightarrow{f} H \xleftarrow{f'} G'$ is defined to be 
\begin{equation}\label{homology operation for a sewing diagram}
\mu_T=f'_*\circ f_!: H_*\bigl(\textup{Map}_*(G,M)\bigr) \rightarrow  H_*\bigl(\textup{Map}_*(H,M)\bigr)
\rightarrow H_*\bigl(\textup{Map}_*(G',M)\bigr),
\end{equation}
which is an operation determined up to sign. This is the string operation associated to the sewing diagram $T$. 

For a labeled sewing diagram $T: G \xrightarrow{f} H \xleftarrow{f'} G'$ and an elementary sewing diagram $T': G' \xrightarrow{g} H' \xleftarrow{g'} G''$ of type (I) or type (II), suppose they are composable and we consider their composition $T''=T'\circ T: G \rightarrow H'' \leftarrow G''$ given in \eqref{composition of sewing diagrams}. Applying the functor $\textup{Map}_*\bigl((\ \cdot\ ), M\bigr)$, we get the following diagram of fibrations over a diagram of their base manifolds: 
\begin{equation*}
\xymatrix{ 
\textup{Map}_*(G,M) & \textup{Map}_*(H,M) \ar[l]_f \ar[r]^{f'}
& \textup{Map}_*(G', M)  \\
& \textup{Map}_*(H'',M) \ar[ul]^{f\circ g''}  \ar[u]_{g''} \ar[r]^{f''}  \ar[dr]_{g'\circ f''}
& \textup{Map}_*(H',M) \ar[u]_g \ar[d]^{g'}  \\
& & \textup{Map}_*(G'',M) }
\end{equation*}
The associated homology diagram with transfers is given by 
\begin{equation*}
\xymatrix{ 
H_*\bigl(\textup{Map}_*(G,M)\bigr) \ar[r]^{f_!} \ar[dr]_{(g'')_!\circ f_!} & H_*\bigl(\textup{Map}_*(H,M)\bigr)  \ar[r]^{f'_*} \ar[d]^{(g'')_!}
& H_*\bigl(\textup{Map}_*(G', M)\bigr) \ar[d]^{g_!} \\
& H_*\bigl(\textup{Map}_*(H'',M)\bigr)    \ar[r]^{f''_*}  \ar[dr]_{(g')_*\circ (f'')_*}
& H_*\bigl(\textup{Map}_*(H',M)\bigr)  \ar[d]^{g'_*}  \\
& & H_*\bigl(\textup{Map}_*(G'',M)\bigr) }
\end{equation*}
Since both triangles and the square commute up to sign, the entire diagram commutes up to sign. Thus a composition of sewing diagrams gives rise to a composition of string operations:
\begin{equation*}
\mu_{T''}=\mu_{T'}\circ\mu_{T}: H_*\bigl(\textup{Map}_*(G,M)\bigr)
\longrightarrow H_*\bigl(\textup{Map}_*(G'',M)\bigr),
\end{equation*}
where the operation $\mu_{T'}$ is an open string product or a coproduct. 

\begin{prop}\label{string operations for sewing diagrams} Let $T: G \xrightarrow{f} H \xleftarrow{f'} G'$ be a sewing diagram labeled by a family of oriented closed submanifolds. For each choice of decomposition 
\begin{equation*}
G \xrightarrow{f_1} H_1 \xrightarrow{f_2} H_2\xrightarrow{f_3} \dotsb \xrightarrow{f_r} H_r=H 
\end{equation*}
of the sewing morphism $f:G \rightarrow H$ into type \textup{(I), (II), (III)} maps, let $T_1, T_2,\dotsc, T_r$ be the elementary sewing diagrams of type \textup{(I)} or \textup{(II)} so that we have $T=T_r\circ T_{r-1}\circ\dotsb\circ  T_1$. Then their corresponding string operations satisfy 
\begin{equation*}
\mu_T=\mu_{T_r}\circ \mu_{T_{r-1}}\circ \dotsb\circ \mu_{T_1}.
\end{equation*}
Thus for any decomposition of $T$ into elementary sewing diagrams, the composition of their corresponding open string products and coproducts is independent of the decomposition of $T$ and is equal to $\mu_T$ given in \eqref{homology operation for a sewing diagram}. 
\end{prop} 

Since a type (I) sewing diagram corresponds to an open string product cobordism and type (II) sewing diagram corresponds to an open string coproduct cobordism, a decomposition of a sewing diagram above corresponds to a decomposition of the corresponding open-closed cobordism $\Sigma$ into an ordered sequence of half-pair-of-pants. We refer to such ordered decomposition into half-pair-of-pants as belonging to $T$. In terms of half-pair-of-pants decomposition, Proposition \ref{string operations for sewing diagrams} can be restated as follows. 

\begin{cor} Let $\Sigma$ be an open-closed cobordism homeomorphic to a disc. Then the associated string operation $\mu_{\Sigma}$ depends only on the sewing diagram $T$, and independent of decompositions of $\Sigma$ into an ordered sequence of half-pair-of-pants belonging to $T$. 
\end{cor}

In the next section, we show that $\mu_\Sigma$ is independent of the choice of the sewing diagram $T$.

\section{Saddle interaction diagrams and open string interactions}

For further discussion, we need to allow open strings to interact at
their internal points. Thus, we introduce more general interaction
diagrams called saddle interaction diagrams. We can then show that
string operations for disc cobordisms are independent of their
decomposition into half-pair-of-pants, and thus independent of sewing
diagrams used for a given disc cobordism. This method also allows us to prove in \S 5 that, when string operations involve open string coproducts, their values are in homology classes of constant open strings.

\subsection{Saddle interaction diagrams and associated string operations}

Let $G$ be a graph consisting of $r$ oriented edges
$e_1,e_2,\dotsc,e_r$. Now suppose that these edges are intersecting at
a single point $p$. Call this configuration $H$. The point $p$ may be
an internal point of an edge, but it can be an end vertex of an
edge. For each edge $e_j$, the point $p$ splits it into an incoming
half edge $e_j'$ and an outgoing half edge $e_j''$. When the point $p$
is an end vertex of $e_j$, then one of the half edges is the
entire edge and the other is a single vertex. To specify how
these $r$ edges interact at $p$, we consider a bijection which we call
a half-edge permutation at $p$:
\begin{equation*}
\sigma_p: \{\textup{ incoming half edges at $p$ }\} \xrightarrow{\cong} \{\textup{ outgoing half edges at $p$ }\}
\end{equation*}
satisfying a condition that $\sigma_p(e_j')\not=e_j''$ for all edges $e_j$. This condition ensures that none of the edges are redundant in the interaction. Given a half-edge permutation $\sigma_p$ for the configuration $H$ above, the set $G'$ of outgoing edges are obtained by joining an incoming half-edge $e_j'$ and an outgoing half-edge $\sigma_p(e_j')$ for each $j$. We have a diagram of the form $G \rightarrow H \leftarrow G'$. Here $G'$ may have a component consisting of a single vertex. This is an example of a saddle interaction diagram, which we define below. 

\begin{defn} (1) Let $G$ be a disjoint union of oriented edges. A \emph{saddle interaction morphism} is a continuous map $f:G \rightarrow H$ onto a forest-like configuration $H$ without loops obtained by intersecting edges of $G$, equipped with a half-edge permutation at each intersection point of $H$. Vertices of $H$ consist of images of vertices of $G$. (The  intersection points of $H$ are not necessarily vertices in our sense. Thus, $H$ is not technically a graph in the usual sense.) 

(2) Such a morphism can be canonically extended to a \emph{saddle interaction diagram} $G \xrightarrow{f} H \xleftarrow{f'} G'$, where $G'$ is obtained from $H$ by recombining half-edges at intersection points according to half-edge permutations. 

(3) Let $\{K_{\lambda}\}_{\lambda\in\Lambda}$ be a family of closed
    oriented submanifolds which are closed under transversal
    intersections. For the empty label $\emptyset$, we set
    $K_{\emptyset}=M$. Suppose vertices of $G$ are labeled by
    $\ell:V(G) \rightarrow \Lambda$. For a vertex $v$ in $H$, if
    $f^{-1}(v)$ contains vertices $v_1,v_2,\dotsc, v_s$ of $G$, then
    the vertex $v$ in $H$ is labeled by the transversal intersection
    $\bigcap_{i=1}^sK_{\ell(v_i)}$. For a vertex $v'$ in $G'$, its
    label is given by the label of the vertex $f'(v')$ in $H$. The
    resulting diagram is a \emph{labeled saddle interaction diagram}.

(4) A \emph{simple saddle interaction diagram} is a saddle interaction
    diagram in which each intersection point of $H$ comes from the 
    intersection of exactly two edges of $G$ at their internal points.
\end{defn}  

Here is an example of a labeled saddle interaction diagram. 

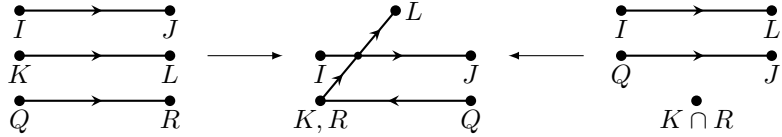
\begin{center}
\begin{tikzpicture}
\path (0,0.6) coordinate (I1) node[below] {$I$};
\fill (I1) circle (2pt);
\path (2,0.6) coordinate (J1) node[below] {$J$};
\fill (J1) circle (2pt);
\path (0,0) coordinate (K1) node[below] {$K$};
\fill (K1) circle (2pt);
\path (2,0) coordinate (L1) node[below] {$L$};
\fill (L1) circle (2pt);
\path (0,-0.6) coordinate (Q1) node[below] {$Q$};
\fill (Q1) circle (2pt);
\path (2,-0.6) coordinate (R1) node[below] {$R$};
\fill (R1) circle (2pt);
\path (4,0) coordinate (I2) node[below] {$I$};
\fill (I2) circle (2pt);
\fill (4.5,0) circle (1.5pt);
\path (6,0) coordinate (J2) node[below] {$J$};
\fill (J2) circle (2pt);
\path (4,-0.6) coordinate (K2) node[below] {$K,R$};
\fill (K2) circle (2pt);
\path (5,0.6) coordinate (L2) node[right] {$L$};
\fill (L2) circle (2pt);
\path (6,-0.6) coordinate (Q2) node[below] {$Q$};
\fill (Q2) circle (2pt);
\path (8,0.6) coordinate (I3)  node[below] {$I$};
\fill (I3) circle (2pt);
\path (10,0.6) coordinate (L3) node[below] {$L$};
\fill (L3) circle (2pt);
\path (8,0) coordinate (Q3) node[below] {$Q$};
\fill (Q3) circle (2pt);
\path (10,0) coordinate (J3)  node[below] {$J$};
\fill (J3) circle (2pt);
\path (9,-0.6) coordinate (KR) node[below] {$K\cap R$};
\fill (KR) circle (2pt);

\path (2.5,0) coordinate (a);
\path (3.5,0) coordinate (b);
\path (6.5,0) coordinate (c);
\path (7.5,0) coordinate (d);

\draw[thick] (I1) -- (J1) (K1) -- (L1) (Q1) -- (R1) 
(I2) -- (J2) (K2) -- (L2) (Q2) -- (K2) 
 (I3) -- (L3)  (Q3) -- (J3);
\draw[->,>=stealth, thick] (1,0) -- ++(0.1,0);
\draw[->,>=stealth, thick] (1,0.6) -- ++(0.1,0);
\draw[->,>=stealth, thick] (1,-0.6) -- ++(0.1,0);
\draw[->,>=stealth, thick] (5,0) -- ++(0.1,0);
\draw[->,>=stealth, thick] (5,-0.6) -- ++(-0.1,0);
\draw[->,>=stealth, thick] (4.25,-0.3) -- ++(0.05,0.05);
\draw[->,>=stealth, thick] (4.75,0.3) -- ++(0.05,0.05);
\draw[->,>=stealth, thick] (9,0) -- ++(0.1,0);
\draw[->,>=stealth, thick] (9,0.6) -- ++(0.1,0);

\draw[-latex] (a) -- (b);  
\draw[-latex]  (d) -- (c);
\end{tikzpicture}
\end{center}

Observe that when exactly two edges intersect, a half-edge permutation
at the intersection point is canonically determined in view of its nontrivial permutation condition. Saddle interaction diagrams are introduced to deal with multiple saddle interactions at internal points involving many open strings. 

When two incoming open strings interact at internal points, recombining half-edges, and become two outgoing open strings, this interaction traces a surface which looks like a saddle. Hence the name saddle interaction in the above definition. 

\begin{rem}\label{difference between diagrams} Sewing diagrams and saddle interaction diagrams are different. They use different labeling scheme, and given $G\rightarrow H$, the method of construction of $G'$ are different, and saddle interaction morphism $f:G\to H$ does not introduce new vertices in $H$, as in type (II) sewing morphisms. The difference becomes most obvious when two open strings interact at their end points, as the following diagrams show. 
\begin{enumerate}
\item[(i)] A labeled sewing diagram with two incoming open strings interacting at their end points carrying the same label $J$. This is the open string product and the corresponding string operation decreases the homological degree by $\dim J$. 

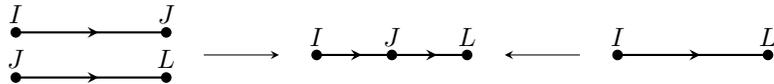
\begin{center}
\begin{tikzpicture}
\path (0,-0.3) coordinate (J0) node[above] {$J$};
\path (2,-0.3) coordinate (L0) node[above] {$L$};
\path (0,0.3) coordinate (I1) node[above] {$I$};
\path (2,0.3) coordinate (J1) node[above] {$J$};
\path (4,0) coordinate (I2) node[above] {$I$};
\path (5,0) coordinate (J2) node[above] {$J$};
\path (6,0) coordinate (L2) node[above] {$L$};
\path (8,0) coordinate (I3) node[above] {$I$};
\path (10,0) coordinate (L3) node[above] {$L$};
\fill (J0) circle (2pt);
\fill (L0) circle (2pt);
\fill (I1) circle (2pt);
\fill (J1) circle (2pt);
\fill (I2) circle (2pt);
\fill (J2) circle (2pt);
\fill (L2) circle (2pt);
\fill (I3) circle (2pt);
\fill (L3) circle (2pt);

\draw[->,>=stealth] (2.5,0) -- (3.5,0);
\draw[->,>=stealth] (7.5,0) -- (6.5,0);
\draw[->,>=stealth,thick] (1,-0.3) -- ++(0.1,0);
\draw[->,>=stealth,thick] (1,0.3) -- ++(0.1,0);
\draw[->,>=stealth,thick] (4.5,0) -- ++(0.1,0);
\draw[->,>=stealth,thick] (5.5,0) -- ++(0.1,0);
\draw[->,>=stealth,thick] (9,0) -- ++(0.1,0);

\draw[thick] (J0) -- (L0) (I1) -- (J1) (I2) -- (L2) (I3) -- (L3);
\end{tikzpicture}
\end{center}

\item[(ii)] A labeled saddle interaction diagram with two incoming open strings interacting at their end points carrying possibly different labels $J$ and $K$. Note that one of the resulting outgoing string is a constant string moving in the transverse intersection $J\cap K$. The associated string operation decreases the homological degree by $\dim M$. 

\begin{center}
\begin{tikzpicture}
\path (0,-0.3) coordinate (K0) node[above] {$K$};
\path (2,-0.3) coordinate (L0) node[above] {$L$};
\path (0,0.3) coordinate (I1) node[above] {$I$};
\path (2,0.3) coordinate (J1) node[above] {$J$};
\path (4,0) coordinate (I2) node[above] {$I$};
\path (5,0) coordinate (J2) node[above] {$J\cap K$};
\path (6,0) coordinate (L2) node[above] {$L$};
\path (8,0.3) coordinate (I3) node[above] {$I$};
\path (10,0.3) coordinate (L3) node[above] {$L$};
\path (9,-0.3) coordinate (K3) node[above] {$J\cap K$};
\fill (K0) circle (2pt);
\fill (L0) circle (2pt);
\fill (I1) circle (2pt);
\fill (J1) circle (2pt);
\fill (I2) circle (2pt);
\fill (J2) circle (2pt);
\fill (L2) circle (2pt);
\fill (I3) circle (2pt);
\fill (L3) circle (2pt);
\fill (K3) circle (2pt);

\draw[->,>=stealth] (2.5,0) -- (3.5,0);
\draw[->,>=stealth] (7.5,0) -- (6.5,0);
\draw[->,>=stealth,thick] (1,-0.3) -- ++(0.1,0);
\draw[->,>=stealth,thick] (1,0.3) -- ++(0.1,0);
\draw[->,>=stealth,thick] (4.5,0) -- ++(0.1,0);
\draw[->,>=stealth,thick] (5.5,0) -- ++(0.1,0);
\draw[->,>=stealth,thick] (9,0.3) -- ++(0.1,0);

\draw[thick] (K0) -- (L0) (I1) -- (J1) (I2) -- (L2) (I3) -- (L3);
\end{tikzpicture}
\end{center}
When $J=K$ in case (ii), $J\cap J$ means the transverse intersection of $J$ with itself. 
\end{enumerate}
However, sewing diagrams without degree $2$ vertices are examples of saddle interaction diagrams. 
\end{rem}

As with labeled sewing diagrams, we can associate a homology operation to labeled saddle interaction diagrams. As before, let $\{K_\lambda\}_{\lambda\in\Lambda}$ be a family of closed oriented submanifolds of $M$. Let $T: G\xrightarrow{f} H \xleftarrow{f'} G'$ be a labeled saddle interaction diagram with the labeling map $\ell: V(G) \rightarrow \Lambda$ for the vertex set $V(G)$. Let $\Delta$ be the set of intersection points of $H$. By taking restricted mapping spaces $\text{Map}_*$ in which vertices are mapped into corresponding submanifolds given by vertex labels, we have the following pull-back squares of fibrations. 
\begin{equation}\label{fibration diagram}
\begin{CD}
\textup{Map}_*(G,M) @<{f}<< \textup{Map}_*(H, M) @>{f'}>> 
\textup{Map}_*(G',M)  \\
@VVV  @VVV  @VVV \\
\textup{Map}_*(f^{-1}(\Delta),M) @<{\bar{f}}<< \textup{Map}_*(\Delta,M)  @>{\bar{f'}}>> \textup{Map}_*((f')^{-1}(\Delta), M) 
\end{CD}
\end{equation}
Here we use the same notation $f$ and $f'$ for induced maps on mapping spaces. In the bottom line, non-vertex points in $\Delta$ which do not carry labels can map anywhere in $M$. The same applies to restricted mapping spaces from sets $f^{-1}(\Delta)$ and $(f')^{-1}(\Delta)$ into $M$. Thus for $p\in \Delta$, if $f^{-1}(p)=\{p_1,p_2,\dotsc,p_k\}\cup\{v_1,v_2,\dotsc,v_r\}$, where $v_i$'s are vertices of $G$ and $p_j$'s are internal points of $G$, then the point $p$ carries the label $\bigcap_{i=1}^rK_{\ell(v_i)}$ and the map $\bar{f}$ is isomorphic to a product over $p\in \Delta$ of the following type of maps of codimension $d(k+r-1)$ for each $p\in\Delta$, where $d=\dim M$:
\begin{equation*}
M^k\times \prod_{i=1}^rK_{\ell(v_i)} \longleftarrow \bigcap_{i=1}^rK_{\ell(v_i)}.
\end{equation*}
Although all labels represent oriented closed submanifolds, the normal bundle to the embedding $\bar{f}$ is only orientable because we do not have preferred ordering of points in $f^{-1}(p)$ for each $p\in \Delta$. Thus all the manifolds in the bottom line in the above diagram are only orientable without preferred orientation. This still gives us a Thom class of the normal bundle to $\bar{f}$ well-defined up to sign. By pulling back the Thom class to the fibration, we can define a transfer map 
\begin{equation*}
f_!: H_*\bigl(\textup{Map}_*(G,M)\bigr) \longrightarrow 
H_*\bigl(\text{Map}(H,M)\bigr)
\end{equation*}
unique up to sign. We then define the string operation associated to labeled saddle interaction diagram $T: G\xrightarrow{f} H \xleftarrow{f'} G'$ by 
\begin{equation*}
\mu_T=f'_*\circ f_!: H_*\bigl(\text{Map}_*(G,M)\bigr) \longrightarrow H_*\bigl(\text{Map}_*(G',M)\bigr).
\end{equation*}
In \eqref{fibration diagram}, we used only intersection points $\Delta$ of $H$.  
We would have obtained the same string operation, up to sign, even if we used $\Delta\cup V(H)$ instead of $\Delta$ in the fibration diagram \eqref{fibration diagram} of restricted mapping spaces, where $V(H)$ is the vertex set of $H$. When $T$ is a sewing diagram without degree $2$ vertices, the associated string operation $\mu_T$ defined earlier coincides with the present string operation for saddle interaction diagrams when $T$ is regarded as a saddle interaction diagram. 

\subsection{Homotopy of saddle interaction diagrams} 

Let $T: G\xrightarrow{f} H \xleftarrow{f'} G'$ be a simple saddle interaction diagram so that each intersection point of $H$ comes from internal points of exactly two edges in $G$. Let $\mu_T$ be the associated string operation. We consider a continuous deformation $\{T_t\}$ of simple saddle interaction diagrams. When the homotopy is through simple saddle interaction diagrams, string operations do not change by homotopy invariance of transfer maps. Suppose this continuous deformation approaches to a multiple saddle interaction diagram or other degenerate diagrams in the "`boundary of the moduli space"' of such diagrams. In this case, when certain transversality conditions are satisfied at the end of the deformation, string operations remain the same during the deformation process. This condition is automatically satisfied for the three deformation processes we consider below as long as the family of closed oriented submanifolds we use as labels of vertices are mutually transversal, which is assumed to be the case. 

We observe that whatever the deformation process is from a simple saddle interaction diagram to a degenerate saddle interaction diagram $T_0: G \xrightarrow{f_0} H_0 \xleftarrow{f'_0} G'$, at each intersection point of $H_0$ the half-edge permutation is canonically determined from simple saddle interaction diagrams by this deformation process. 

We describe three types of deformation processes below. Let $T: G\xrightarrow{f} H \xleftarrow{f'} G'$ be a simple saddle interaction diagram with connected $H$. Let $\Delta\subset H$ be the set of intersection points, which come from internal points of edges in $G$ via $f$. Here we note that when the set $\Delta$ of intersection points consists of $r$ points, then $G$ has $r+1$ disjoint edges since $H$ is connected and has no loops. The diagram relevant for string operation $\mu_T$ is the following: 
\begin{equation}\label{simple saddle diagram}
\begin{CD}
\textup{Map}_*(G,M) @<{f}<< \textup{Map}_*(H,M) @>{f'}>> \textup{Map}_*(G',M) \\
@VVV @VVV @.  \\
\textup{Map}_*(f^{-1}(\Delta),M) @<{\bar{f}}<< \textup{Map}_*(\Delta, M) @. \\
@A{\cong}AA  @A{\cong}AA @. \\
(M\times M)^r @<{\phi\times\dotsb\times\phi}<< M^r @. 
\end{CD}
\end{equation}
where the downward vertical maps are induced by restrictions. 

(1) (Deformation to a sewing diagram)  Let $\{T_t: G \xrightarrow{f_t} H_t \xleftarrow{f'_t} G'_t\}$ be a homotopy between simple saddle interaction diagrams $T_t$ for $0<t\le 1$ and a sewing diagram $T_0$ with only degree $1$ and $3$ vertices at $t=0$ described as follows.  In this homotopy at each intersection point of $H_1$, one of the two edges slides and the intersection point coincides at $t=0$ with one of its end vertices. Let $\Delta_0=\{v_1,v_2,\dotsc,v_r\}$ be degree $3$ vertices in $H_0$ carrying labels $K_1,K_2,\dotsc,K_r$. They are type (III) vertices. Then for each $1\le j\le r$, the inverse image $f_0^{-1}(v_j)$ consists of an internal point of an edge and an end vertex of another edge carrying the label $K_j$. Thus a diagram for the string operation $\mu_{T_0}$ is 
\begin{equation*}
\begin{CD}
\textup{Map}_*(G,M) @<{f_0}<< \textup{Map}_*(H_0,M) @>{f'_0}>> \textup{Map}_*(G_0',M) \\
@VVV @VVV @.  \\
\textup{Map}_*(f_0^{-1}(\Delta_0),M) @<{\bar{f}_0}<< \textup{Map}_*(\Delta_0, M) @. \\
@A{\cong}AA  @A{\cong}AA @. \\
\prod_{j=1}^r(K_j\times M) @<{\prod_j(1\times i_j)\phi_j}<< 
\prod_{j=1}^rK_j  @. 
\end{CD}
\end{equation*}
where for $1\le j\le r$, the map $\phi_j:K_j \rightarrow K_j\times K_j$ is a diagonal map and $\iota_j:K_j \rightarrow M$ is an inclusion map. We compare the map between base manifolds of the above diagram with the one for simple saddle interaction. 
\begin{equation*}
\begin{CD}
(M\times M)\times \cdots \times(M\times M) @<{\ \ \ \ \ \ \ 
 \phi\times \cdots\times\phi\ \ \ \ \ \ \ }<< M\times \cdots\times M=M^r \\
@AA{(\iota_1\times 1)\times \cdots\times(\iota_r\times 1)}A 
@AA{\iota_1\times\cdots\times\iota_r}A  \\
(K_1\times M)\times\cdots\times(K_r\times M) 
@<{(1\times \iota_1)\phi_1\times\cdots\times(1\times\iota_r)\phi_r}<<
K_1\times\cdots\times K_r 
\end{CD}
\end{equation*}
Here, the top line is the inclusion map of base manifolds for the fibration in \eqref{simple saddle diagram}. 
Since this is a pull-back diagram, the normal bundle of the embedding $\phi\times \cdots\times \phi$ in the top row restricts to the normal bundle of the embedding in the bottom horizontal map. Thus, Thom classes for these normal bundles coincide up to sign under the restriction map $\iota_1\times\dotsb\times\iota_r$. The sign ambiguity is present because an ordering of vertices $v_1,\dotsc, v_r$ is not specified. Consequently the orientations of orientable manifolds $M^r$ and $\prod^r_{i=1}K_i$ are not specified. This means that the string operations associated to the deformation of saddle interaction diagrams does not change. In particular, a two-string saddle operation coincides with an open string coproduct followed by an open string product. 

(2) (Deformation to a simultaneous internal point interaction) 
In this deformation process from a simple saddle interaction diagram $T_1:G \to H_1 \gets G_1'$ with connected $H_1$, all the internal intersection points in $H_1$ at $t=1$ come together to the same internal interaction point $p$ in $H_0$ at $t=0$. Since $H_t$ is connected and without loops, if its set of interaction points $\Delta_t$ consists of $r$ points, then $G$ has $r+1$ edges and $f_0^{-1}(p)$ consists of $r+1$ internal points of $G$, one for each edge.  The string operation $\mu_{T_0}$ comes from the following diagram: 
\begin{equation*}
\begin{CD}
\textup{Map}_*(G,M) @<{f_0}<< \textup{Map}_*(H_0,M) @>{f'_0}>> \textup{Map}_*(G_0',M) \\
@VVV @VVV @.  \\
\textup{Map}_*(f_0^{-1}(\Delta_0),M) @<{\bar{f}_0}<< \textup{Map}_*(\Delta_0, M) @. \\
@A{\cong}AA  @A{\cong}AA @. \\
M^{r+1} @<{\phi_{r+1}}<< M @. 
\end{CD}
\end{equation*}
Here $\phi_{r+1}$ is a diagonal map to a product of $r+1$ copies of $M$. Comparing the base manifold embedding of the simple saddle interaction diagram \eqref{simple saddle diagram} with the above base manifold embedding, we have the following pull-back diagram: 
\begin{equation*}
\begin{CD}
(M\times M)^r 
@<{\phi\times\cdots\times\phi}<< 
M^r \\
@AAA   @AA{\phi_r}A  \\
M^{r+1}  @<{\phi_{r+1}}<<  M 
\end{CD}
\end{equation*}
Here the left vertical map depends on the individual degeneration process. Again the normal bundle of the embedding $\phi\times\cdots\times \phi$ in the top row restricts to the normal bundle of the bottom embedding $\phi_{r+1}$, and consequently, the Thom classes correspond under the restriction map, up to sign. Type (1) and (2) deformations are used to prove Theorem D (which is Theorem \ref{independence of decomposition} below). 

(3) (Deformation to a simultaneous vertex interaction) In this process, all edges in $H_1$ slide toward one of their end vertices, and at $t=0$, all of these end vertices coincide and become the only intersection point in $H_0$. Let $\Delta_0=\{v\}$, and suppose vertices in the inverse image $f_0^{-1}(v)=\{v_1,v_2,\dotsc,v_{r+1}\}$ carry labels $K_1,K_2,\dotsc,K_{r+1}$. Then, the simultaneous intersection vertex $v$ carries the label $K_1\cap K_2\cap\cdots\cap K_{r+1}$, where the intersection is taken transversally (when $K_i=K_j$ for some $i\ne j$, we deform $K_i$ to obtain the transversal intersection of $K_i$ with itself). In this case, the diagram computing the string operation $\mu_{T_0}$ is the following one.  
\begin{equation*}
\begin{CD}
\textup{Map}_*(G,M) @<{f_0}<< \textup{Map}_*(H_0,M) @>{f'_0}>> \textup{Map}_*(G',M) \\
@VVV @VVV @.  \\
\textup{Map}_*(f_0^{-1}(\Delta_0),M) @<{\bar{f}_0}<< \textup{Map}_*(\Delta_0, M) @. \\
@A{\cong}AA  @A{\cong}AA @. \\
K_1\times\cdots\times K_{r+1}  @<{\phi_{r+1}}<< 
K_1\cap \cdots\cap K_{r+1}
\end{CD}
\end{equation*}
To see that Thom classes behave well when degeneration takes place at $t=0$, we check that the following diagram of inclusion maps of base manifolds of fibrations is a pull-back diagram, which is obvious: 
\begin{equation*}
\begin{CD}
(M\times M)^r 
@<{\phi\times\cdots\times\phi}<< 
M^r \\
@AAA   @AA{\phi_r\circ \iota}A  \\
K_1\times K_2\times \cdots\times K_{r+1} 
@<{\phi_{r+1}}<<  K_1\cap K_2\cap\cdots\cap K_{r+1},
\end{CD}
\end{equation*}
where $\iota:K_1\cap\cdots\cap K_{r+1} \rightarrow M$ is the inclusion map. 
The left vertical map depends on the individual deformation process of edge configuration. Type (3) deformation is used in \S5 to extract qualitative and essential features of general saddle interactions described in Theorem C. 

In the above three deformation processes, the family of tubular neighborhoods of the embedding 
$\textup{Map}_*(G,M) \xleftarrow{f_t} \textup{Map}_*(H_t,M)$ for $0\le t\le 1$ fit together, so the Thom class for the individual $t$ comes from the same Thom class for the family. Hence by the homotopy invariance of the transfer, string operations remain the same through the deformation process. 

As an application of these deformation processes, we can show that the string operation $\mu_{\Sigma}$ for an open-closed disc cobordism is independent of ordered decompositions of $\Sigma$ into half-pair-of-pants. 
Recall that in Proposition \ref{string operations for sewing diagrams}, we showed that the string operation $\mu_{\Sigma}$ depends only on the sewing diagram $T$ used for decomposition of $\Sigma$. We now show that $\mu_{\Sigma}$ is independent of sewing diagrams and depends only on $\Sigma$. 

\begin{thm}\label{independence of decomposition} Let $\Sigma$ be an open cobordism homeomorphic to a disc. Then the string operation $\mu_{\Sigma}$ associated to $\Sigma$, determined up to sign, is  independent of its decomposition into half-pair-of-pants.  
\end{thm}
\begin{proof}  Suppose we have two decompositions of $\Sigma$ into half-pair-of-pants belonging to sewing diagrams $T$ and $T'$. By Proposition \ref{string operations for sewing diagrams}, the string operations associated to these decompositions of $\Sigma$ depend only on associated sewing diagrams $T$ and $T'$, up to sign. We consider the canonical decompositions of sewing diagrams $T,T'$ given in Proposition \ref{canonical decomposition of sewing diagrams}: 
\begin{equation*}
T=T_2\circ T_3\circ T_1,\qquad 
T'=T_2'\circ T_3' \circ T_1'.
\end{equation*}
Since $T_1$ and $T_1'$ have sewing graphs consisting of edges with only type (I) degree $2$ vertices, and they correspond to combining incoming open strings which are next to each other along the boundary of $\Sigma$, they depend only on $\Sigma$, and hence $T_1=T_1'$. Similarly, $T_2$ and $T_2'$ have sewing graphs consisting of edges with only type (II) degree $2$ vertices, and represent processes of splitting outgoing open strings. Thus they depend only on adjacent outgoing open strings in $\Sigma$. Hence $T_2=T_2'$. Since $\mu_{T}=\mu_{T_2}\circ\mu_{T_3}\circ\mu_{T_1}$ and $\mu_{T'}=\mu_{T_2'}\circ\mu_{T_2'}\circ\mu_{T_1'}$, we only have to show $\mu_{T_3}=\mu_{T_3'}$.  

To see this, let $H_3$ and $H_3'$ be sewing graphs for the sewing diagrams $T_3$ and $T_3'$. So $H_3$ and $H_3'$ are forests with only  vertices of degree $1$ or $3$. Sewing diagrams $T_3$ and $T_3'$ represent the same open cobordism obtained from $\Sigma$ by combining adjacent incoming open strings and adjacent outgoing open strings. Thus, $T_3$ and $T_3'$ are of the form $T_3: G \rightarrow H_3 \leftarrow G'$ and $T_3': G \rightarrow H_3' \leftarrow G'$ for the same disjoint union of labeled edges $G$ and $G'$. 

Now we can continuously deform the sewing diagram $T_3$ through simple saddle interaction diagrams (backward deformation of type (1) above) to a simultaneous saddle interaction diagram $\widetilde{T}_3$ with a single internal multiple intersection point (deformation of type (2) above). Similarly, starting with a sewing diagram $T_3'$, we can deform it through simple saddle interaction diagrams to a simultaneous saddle interaction diagram $\widetilde{T}_3'$  (backward deformation of type (1) again followed by type (3) deformation). As we saw earlier, the associated string operations remain the same throughout these deformations of saddle interaction diagrams. Since the saddle interaction diagrams $\widetilde{T}_3$ and $\widetilde{T}_3'$ are homeomorphic, they induce the same string operation and we have $\mu_{T_3}=\mu_{\widetilde{T}_3}=\mu_{\widetilde{T}_3'}=\mu_{T_3'}$. Hence string operations for $T_3$ and $T_3'$ must coincide. 
This completes the proof. 
\end{proof} 

Figure 7 describes an example of a deformation between two sewing graphs $H_3$ and $H_3'$ without degree $2$ vertices through a saddle interaction graph $\widetilde{H}_3$ with one multiple intersection point $p$. Each of the three saddle interaction graphs below describes an interaction of three incoming open strings with labels $(I,J), (K,L), (Q,R)$, producing three outgoing open strings with labels $(I,L), (K,R), (Q,J)$. During the deformation on the left, edges $KL$ and $QR$ slides up or down and two intersection points coincide at $p$. Similarly for the deformation on the right. At the intersection point $p$, a half edge permutation $\sigma_p$ is given so that we obtain edges $IL$, $KR$, and $QJ$ after interaction. 

In type (1) deformation above, we considered deformations of generic saddle interaction diagrams into sewing diagrams without degree 2 vertices. For general sewing diagrams, those with type (II) degree 2 vertices (which correspond to open string coproducts) are amenable for deformations in which type (II) degree 2 vertices merge into a single multiple type (II) vertex, corresponding to a degenerate iterated open string coproduct. In the next section, we use these deformations to deduce constancy results for (iterated) open string coproducts. 

\begin{figure}
\begin{center}
\begin{tikzpicture}[decoration={snake, amplitude=0.4mm, segment length=2mm, post length=0.1cm}]

\path (0,0) coordinate (I1) node[below] {$I$} ;
\fill (I1) circle (2pt) ; 
\path (2,0) coordinate (J1) node[below] {$J$} ;
\fill (J1) circle (2pt); 
\path (0.5,0) coordinate (K1) node[below] {$K$} ;
\fill (K1) circle (2pt); 
\path (1,1) coordinate (L1) node[right] {$L$} ;
\fill (L1) circle (2pt); 
\path (1.2,0) coordinate (Q1) node[above] {$Q$} ;
\fill (Q1) circle (2pt); 
\path (1.8,-1) coordinate (R1) node[right] {$R$} ;
\fill (R1) circle (2pt); 
\path (3.5,0) coordinate (I2) node[below] {$I$} ;
\fill (I2) circle (2pt); 
\path (5.5,0) coordinate (J2) node[below] {$J$} ;
\fill (J2) circle (2pt); 
\path (4,-1) coordinate (K2) node[left] {$K$} ;
\fill (K2) circle (2pt); 
\path (5,1) coordinate (L2) node[right] {$L$} ;
\fill (L2) circle (2pt); 
\path (4,1) coordinate (Q2) node[left] {$Q$} ;
\fill (Q2) circle (2pt); 
\path (5,-1) coordinate (R2) node[right] {$R$} ;
\fill (R2) circle (2pt); 
\path (4.5,0) coordinate (p) node[below right] {$p$};
\path (7,0) coordinate (Q3) node[left] {$Q$} ;
\fill (Q3) circle (2pt); 
\path (8.5,0) coordinate (R3) node[below] {$R$} ;
\fill (R3) circle (2pt); 
\path (7.5,0) coordinate (K3) node[below] {$K$} ;
\fill (K3) circle (2pt); 
\path (8,1) coordinate (L3) node[right] {$L$} ;
\fill (L3) circle (2pt); 
\path (6.5,1) coordinate (I3) node[left] {$I$} ;
\fill (I3) circle (2pt); 
\path (7.5,-1) coordinate (J3) node[right] {$J$} ;
\fill (J3) circle (2pt); 

\path (2.3,0) coordinate (a) ;
\path (3.2,0) coordinate (b);
\path (5.8,0) coordinate (c) ;
\path (6.5,0) coordinate (d) ;

\draw[-latex, decorate] (a) -- (b) ;
\draw[-latex, decorate] (d) -- (c) ;

\draw[thick] (I1) -- (J1)  (K1) -- (L1)  (Q1) -- (R1)  
(I2) -- (J2)  (Q2) -- (R2)  (K2) -- (L2)  
(Q3) -- (R3)  (I3) -- (J3)  (K3) -- (L3);

\draw[->, >=stealth, thick] (1,0) -- ++(0.1,0);
\draw[->, >=stealth, thick] (0.75,0.5) -- ++(0.05,0.1);
\draw[->, >=stealth, thick] (1.5,-0.5) -- ++(0.06,-0.1);
\draw[->, >=stealth, thick] (4,0) -- ++(0.1,0);
\draw[->, >=stealth, thick] (5,0) -- ++(0.1,0);
\draw[->, >=stealth, thick] (4.25,0.5) -- ++(0.05,-0.1);
\draw[->, >=stealth, thick] (4.75,-0.5) -- ++(0.05,-0.1);
\draw[->, >=stealth, thick] (4.25,-0.5) -- ++(0.05,0.1);
\draw[->, >=stealth, thick] (4.75,0.5) -- ++(0.05,0.1);
\draw[->, >=stealth, thick] (8,0) -- ++(0.1,0);
\draw[->, >=stealth, thick] (6.75,0.5) -- ++(0.05,-0.1);
\draw[->, >=stealth, thick] (7.25,-0.5) -- ++(0.05,-0.1);
\draw[->, >=stealth, thick] (7.75,0.5) -- ++(0.05,0.1);

\path node at (1,-1.5) {$H_3$};
\path node at (4.5,-1.5) {$\widetilde{H}_3$};
\path node at (7.5,-1.5) {$H_3'$};

\draw node[text width=11cm] at (4.5, -3) 
{\textsc{Figure 7.} A deformation between sewing graphs $H_3$ and $H_3'$ through saddle interaction graph $\widetilde{H}_3$. The corresponding string operations remain the same throughout this deformation process. This leads to the proof that the disc string operation is independent of its half-pair-of-pants decompositions (Theorem 3.3). };


\end{tikzpicture}
\end{center}
\end{figure}

\section{Open string (co)products and constant homology classes}

Let $I,J,K,L,\dotsc,$ be closed oriented submanifolds of $M$ and assume that they are transversal to each other. The homology class of constant open strings $[I]\in H_*(P_{II})$ is the unit in the algebra $H_*(P_{II})$ with respect to the open string product. For a path space $P_{IJ}$, its homology group $H_*(P_{IJ})$ does not have the structure of an algebra, but the homology class $[I\cap J]\in H_*(P_{IJ})$ of constant open strings plays a very interesting role in the open string coproduct.

\subsection{Properties of the open string coproduct}

Let $\mu_J: H_*(P_{IJ})\otimes H_*(P_{JK}) \rightarrow H_*(P_{IK})$ be the open string product map along $J$. Orient transversal intersections $I\cap J$ and $J\cap K$ in any way, and let $[I\cap J]\in H_*(P_{IJ})$ and $[J\cap K]\in H_*(P_{JK})$ be the homology classes of constant open strings. 

\begin{lem}\label{product with fundamental class} \textup{(1)}  For $a\in H_*(P_{JK})$, 
\begin{equation*}
\mu_J([I\cap J]\otimes a)=\pm(\iota_I)_*(\iota_J)_!(a)\in H_*(P_{IK}),
\end{equation*}
where $\iota_I$ and $\iota_J$ are obvious inclusion maps $P_{JK}\xleftarrow{\iota_J} P_{I\cap J, K} 
\xrightarrow{\iota_I} P_{IK}$.

\textup{(2)} For $b\in H_*(P_{IJ})$, 
\begin{equation*}
\mu_J(b\otimes [J\cap K])=\pm(\iota_K)_*(\iota_J)_!(b)\in H_*(P_{IK}),
\end{equation*}
where $\iota_J$ and $\iota_K$ are inclusion maps 
$P_{IJ} \xleftarrow{\iota_J}  P_{I,J\cap K} \xrightarrow{\iota_K}
P_{IK}$. The sign $\pm$ depends on chosen orientations. 
\end{lem}
\begin{proof} We consider the following homotopy commutative diagram, where the left square strictly commutes. 
\begin{equation*}
\begin{CD}
P_{IJ}\times P_{JK} @<{j_J}<< P_{IJ}\underset{J}{\times}P_{JK}
@>{\iota_J}>> P_{IK}  \\
@AA{s\times 1}A @AA{s'}A  @|  \\
(I\cap J)\times P_{JK} @<{(p_0,\iota_J)}<< P_{I\cap J,K} @>{\iota_I}>> P_{IK}, 
\end{CD}
\end{equation*}
where $s:I\cap J \rightarrow P_{IJ}$ is the inclusion map,  $s'(\gamma)=(c_{\gamma(0)},\iota_J(\gamma))$ for $\gamma\in P_{I\cap J,K}$, and $p_0(\gamma)=\gamma(0)$. Here, $c_x$ is the constant path at $x\in M$. Since the left square is a pull-back diagram, the associated homology diagram with transfers commutes up to sign:
\begin{equation*}
\begin{CD}
H_*(P_{IJ}\times P_{JK}) @>{(j_J)_!}>> H_*(P_{IJ}\underset{J}{\times}P_{JK})
@>{(\iota_J)_*}>> H_*(P_{IK})  \\
@AA{(s\times 1)_*}A @AA{s'_*}A  @|  \\
H_*\bigl((I\cap J)\times P_{JK}\bigr) @>{(p_0,\iota_J)_!}>> 
H_*(P_{I\cap J,K}) @>{(\iota_I)_*}>> H_*(P_{IK}) 
\end{CD}
\end{equation*}
Let $\pi_2: (I\cap J)\times P_{JK} \rightarrow P_{JK}$ be a projection map. 
Since $(p_0,\iota_J)_!([I\cap J]\times a)=\pm (p_0,\iota_J)_!(\pi_2)_!(a)=\pm(\iota_J)_!(a)$ for $a\in H_*(P_{JK})$, the commutativity of the above diagram implies that $\mu_J(s_*([I\cap J])\otimes a)=\pm(\iota_I)_*(\iota_J)_!(a)$. This proves (1). The part (2) can be proved similarly. 
\end{proof}

It turns out that the product with the ``fundamental'' homology class of the constant paths $[I\cap J]\in H_*(P_{IJ})$ has the effect of mapping arbitrary homology classes in $H_*(P_{JK})$ to constant homology classes in $H_*(P_{IK})$ which are in the image from $H_*(I\cap K)$. We now show that this property of the open string product is related to a property of the open string coproduct that the image of the open string coproduct consists of homology classes of constant paths. The idea here is to deform type (II) sewing diagrams. Let $\varphi_I: H_*(P_{JK}) \rightarrow H_*(P_{JI})\otimes H_*(P_{IK})$ be the open string coproduct split along $I$. Let $p_0:P_{IJ} \rightarrow I$ and $p_1:P_{IJ} \rightarrow J$ be projections to initial and end points of open strings. Namely $p_0(\gamma)=\gamma(0)$ and $p_1(\gamma)=\gamma(1)$ for $\gamma\in P_{IJ}$. 

For a geometric meaning of the next theorem, which may be more illuminating than the proof itself, see the discussion after the proof.

\begin{thm}\label{image of open coproducts}  \textup{(1)} The image of the open string coproduct $\varphi_J$ consists of homology classes of constant paths. Namely, the coproduct map $\varphi_J$ factors as follows.
\begin{equation*}
\xymatrix{
H_*(P_{JK}) \ar[r]^{\varphi_J\ \ \ \ \ \ \ } \ar[dr] & H_*(P_{JI})\otimes H_*(P_{IK}) \\
& H_*(J\cap I)\otimes H_*(I\cap K) \ar[u]_{s_*\otimes s_*}
}
\end{equation*}

\textup{(2)}  For $a\in H_*(P_{JK})$, the open string product with $[I\cap J]\in H_*(P_{IJ})$ gives a homology class of constant open strings, as in the following formula and the commutative diagram\textup{:}
\begin{equation*}
\mu_J([I\cap J]\otimes a)=\pm s_*p_{1*}(\iota_K')_!(a), \quad
\begin{CD}
H_*(P_{JK})  @>{\mu_J([I\cap J]\otimes (\cdot))}>> H_*(P_{IK}) \\
@VV{(\iota_K')_!}V   @AA{s_*}A  \\
H_*(P_{J,I\cap K})  @>{p_{1*}}>>  H_*(I\cap K)
\end{CD}
\end{equation*}
Similarly for $b\in H_*(P_{IJ})$, 
\begin{equation*}
\mu_J(b\otimes[J\cap K])=\pm s_*p_{0*}(\iota_K)_!(b),\quad
\begin{CD}
H_*(P_{IJ}) @>{\mu_J((\cdot)\otimes[J\cap K])}>> 
H_*(P_{IK})  \\
@VV{(\iota_I)_!}V   @AA{s_*}A  \\
H_*(P_{I\cap K,J}) @>{p_{0*}}>>  H_*(I\cap K)
\end{CD}
\end{equation*}

\textup{(3)}  We have the following diagram commutative up to sign, 
\begin{equation*}
\begin{CD}
H_*(P_{JK}) @>{[I\cap J]\otimes(\cdot)}>> [I\cap J]\otimes H_*(P_{JK})  \\
@VV{\varphi_I}V    @VV{\mu_J}V  \\
H_*(P_{JI})\otimes H_*(P_{IK}) @>{(\epsilon_{JI})_*\otimes 1}>> H_*(P_{IK}),
\end{CD}
\end{equation*}
where $\epsilon_{JI}: P_{JI} \rightarrow \textup{pt}$ is the constant map, and the top horizontal map sends $a\in H_*(P_{JK})$ to $[I\cap J]\otimes a$. Here the image of vertical maps are constant homology classes by \textup{(1)} and \textup{(2)}. Similarly, we have the following commutative diagram up to sign.
\begin{equation*}
\begin{CD}
H_*(P_{IJ}) @>{(\cdot)\otimes[J\cap K]}>> H_*(P_{IJ})\otimes [J\cap K]  \\
@VV{\varphi_K}V  @VV{\mu_J}V  \\
H_*(P_{IK})\otimes H_*(P_{KJ})  @>{1\otimes (\epsilon_{KJ})_*}>> H_*(P_{IK})
\end{CD}
\end{equation*}
\end{thm}

\begin{figure}
\begin{center}
\begin{tikzpicture}
\path (0,0) node[below] {$J$};
\path (0,0) node[above] {$0$};
\fill (0,0) circle (2pt);
\path (1,0) node[below] {$I$};
\path (1,0) node[above] {$t$};
\fill (1,0) circle (2pt);
\path (2,0) node[below] {$K$};
\path (2,0) node[above] {$1$};
\fill (2,0) circle (2pt);

\path (4,0) node[below] {$I\cap J$};
\path (4,0) node[above] {$t=0$};
\fill (4,0) circle (2pt);
\path (6,0) node[below] {$K$};
\path (6,0) node[above] {$1$};
\fill (6,0) circle (2pt);

\path (8,0) node[below] {$J$};
\path (8,0) node[above] {$0$};
\fill (8,0) circle (2pt);
\path (10,0) node[below] {$I\cap K$};
\path (10,0) node[above] {$t=1$};
\fill (10,0) circle (2pt);

\draw[thick] (0,0) -- (2,0)  (4,0) -- (6,0)  (8,0) -- (10,0);

\draw[->,>=stealth] (0.5,0) -- ++(0.1,0);
\draw[->,>=stealth] (1.5,0) -- ++(0.1,0);
\draw[->,>=stealth] (5,0) -- ++(0.1,0);
\draw[->,>=stealth] (9,0) -- ++(0.1,0);

\path (-0.7,0) node {$(a)$};
\path (3,0) node {$(b)$};
\path (7.3,0) node {$(c)$};

\draw node[text width=11cm] at (5,-2)
{\textsc{Figure 8.} In the sewing diagram for the open string coproduct, the middle interacting vertex can be anywhere along the open string without changing the string operation. It can be at $t=0$ as in $(b)$, or at $t=1$ as in $(c)$. This observation leads to the proof of the constancy of the value of the open string coproduct (Theorem 4.2). };

\end{tikzpicture}
\end{center}
\end{figure}

\begin{proof} (1) We consider splitting of paths in $P_{JK}$ at time $t\in [0,1]$ by $I$, and we set  $P_{JI}\underset{I}{\overset{t}{\times}}P_{IK}$ be the set of paths $\gamma\in P_{JK}$ such that $\gamma(t)\in I$. The corresponding diagram is given in Figure 8 (a) in which the coordinate $t$ of the middle vertex can move freely in the interval $[0,1]$. When $t=0,1$, we have the following situation corresponding to diagrams in Figure 8 (b) and (c). 
\begin{equation*}
\begin{cases}
P_{JI}\underset{I}{\overset{0}{\times}}P_{IK}\cong P_{I\cap J,K}& 
\text{at $t=0$}, \\ P_{JI}\underset{I}{\overset{1}{\times}}P_{IK}\cong P_{J,I\cap K}&
\text{ at $t=1$}. 
\end{cases}
\end{equation*}
We consider the following commutative diagram, where $j_0(\gamma)=(c_{\gamma(0)}, \iota_I(\gamma))$ for $\gamma\in P_{I\cap J,K}$, and $j_1(\gamma)=(\iota_I'(\gamma), c_{\gamma(1)})$ for $\gamma\in P_{J,I\cap K}$. 
\begin{equation*}
\xymatrix{
 & P_{I\cap J,K} \ar[dl]^{\iota_0=\iota_J} \ar[dr]^{j_0} 
\ar[r]^{(p_0,1)\ \ \ \ \ \ \ } 
& (I\cap J)\times P_{I\cap J,K} \ar[d]^{s\times \iota_I} 
\ar[r]^{\ \ \ \ \ \ \ \pi_2}  
& P_{I\cap J,K} \ar[d]^{\iota_I} \\
P_{JK} & P_{JI}\underset{I}{\overset{t}{\times}}P_{IK} 
\ar[l]_{\ \ \ \iota_t\ \ \ \ \ } \ar[r]^{j_t} 
& P_{JI}\times P_{IK}   \ar[r]^{\ \ \pi_2} & P_{IK}  \\
& P_{J,I\cap K} \ar[ul]_{\iota_1=\iota_K'} \ar[ur]^{j_1}
\ar[r]^{(1,p_1)\ \ \ \ \ \ \ }  
& P_{J,I\cap K}\times (I\cap K)  \ar[u]_{\iota_I'\times s} 
\ar[r]^{\ \ \ \ \ \ \ \ \pi_2} & I\cap K \ar[u]_{s}
}
\end{equation*}
By homotopy invariance of the transfer, in the associated homology diagram with transfers we have 
\begin{equation*}
\varphi_J=(j_t)_*(\iota_t)_!=(j_0)_*(\iota_0)_!
=(j_1)_*(\iota_1)_!: H_*(P_{JK}) \rightarrow 
H_*(P_{JI})\otimes H_*(P_{Ik}).
\end{equation*}
Using factorization of $j_0$ and $j_1$ as above, the third column of the diagram implies 
\begin{align*}
\text{Im}\,\varphi_J&\subset s_*\bigl(H_*(I\cap J)\bigr)\otimes (\iota_I)_*\bigl(H_*(P_{I\cap J, K})\bigr) \\
\text{Im}\,\varphi_J&\subset (\iota_I)_*\bigl(H_*(P_{J,I\cap K})\bigr)
\otimes s_*\bigl(H_*(I\cap K)\bigr).
\end{align*}
Hence we have $\text{Im}\,\varphi_J\subset s_*\bigl(H_*(I\cap J)\bigr)\otimes s_*\bigl(H_*(I\cap K)\bigr)$. This proves (1). 

Following the upper perimeter from left to right, we get for $a\in H_*(P_{JK})$,
\begin{equation*}
(\iota_I)_*(\pi_2)_*(p_0,1)_*(\iota_J)_!(a)=(\iota_I)_*(\iota_J)_!(a)
=\pm \mu_J([I\cap J]\otimes a),
\end{equation*}
using $\pi_2\circ(p_0,1)=1$ and Lemma \ref{product with fundamental class}. Following the middle row gives 
\begin{equation*}
(\pi_2)_*(j_t)_*(\iota_t)_!(a)=\bigl((\epsilon_{JI})_*\otimes 1\bigr)
\varphi_J(a).
\end{equation*}
Following the lower perimeter gives 
\begin{equation*}
s_*(\pi_2)_*(1,p_1)_*(\iota_K')_!(a)=s_*p_{1*}(\iota_K')_!(a).
\end{equation*}
Since the associated homology diagram with transfers commutes up to sign, we have 
\begin{equation*}
\mu_J([I\cap J]\otimes a)
=\pm \bigl((\epsilon_{JI})_*\otimes 1\bigr)\varphi_J(a)
=\pm s_*p_{1*}(\iota_K')_!(a).
\end{equation*}
This proves the first parts of (2) and (3). The second parts of (2) and (3) can be proved using the following diagram. 
\begin{equation*}
\xymatrix{
 & P_{I\cap K,J} \ar[dl]^{\iota_0=\iota_J} \ar[dr]^{j_0} 
\ar[r]^{(p_0,1)\ \ \ \ \ \ \ } 
& (I\cap K)\times P_{I\cap K,J} \ar[d]^{s\times \iota_K} 
\ar[r]^{\ \ \ \ \ \ \ \pi_1}  
& I\cap K \ar[d]^{s} \\
P_{IJ} & P_{IK}\underset{K}{\overset{t}{\times}}P_{KJ} 
\ar[l]_{\ \ \ \iota_t\ \ \ \ \ \ } \ar[r]^{j_t} 
& P_{IK}\times P_{KJ}   \ar[r]^{\ \ \pi_1} & P_{IK}  \\
& P_{I,K\cap J} \ar[ul]_{\iota_1=\iota_J'} \ar[ur]^{j_1}
\ar[r]^{(1,p_1)\ \ \ \ \ \ \ }  
& P_{I,K\cap J}\times (K\cap J)  \ar[u]_{\iota_K'\times s} 
\ar[r]^{\ \ \ \ \ \ \ \ \pi_1} & P_{I,K\cap J} \ar[u]_{\iota_K'}
}
\end{equation*}
This completes the proof.
\end{proof}

\begin{example} As an example, we consider the following composition of open string product maps for oriented closed submanifolds $I,J,L,Q$:
\begin{equation*}
\mu_{JL}=\mu_L\circ(\mu_J\otimes 1): 
H_*(P_{IJ})\otimes H_*(P_{JL})\otimes H_*(P_{LQ}) \rightarrow  H_*(P_{IQ}).
\end{equation*}
Then we claim that for homology classes $a\in H_*(P_{IJ})$, $[J\cap K\cap L]\in H_*(P_{JL})$, and $b\in H_*(P_{LQ})$ for an oriented closed manifold $K$, the following open string product is a constant homology class: 
\begin{equation*}
\mu_{JL}(a\otimes[J\cap K\cap L]\otimes b)\in 
s_*\bigl(H_*(I\cap K\cap Q)\bigr). 
\end{equation*}
To see this, we first note that for fundamental homology classes of constant path $[J\cap K]\in H_*(P_{JK})$ and $[K\cap L]\in H_*(P_{KL})$, we have $\mu_K([J\cap K]\otimes [K\cap L])=[J\cap K\cap L]$ with respect to suitably chosen orientations. Here, intersections are taken as transverse intersections. By the associativity of open string products, we have 
\begin{equation*}
\mu_{JL}(a\otimes[J\cap K\cap L]\otimes b)=\pm \mu_K\circ(\mu_J\otimes\mu_L)(a\otimes [J\cap K]\otimes[K\cap L]\otimes b).
\end{equation*}
Since by Proposition \ref{image of open coproducts}, we have 
\begin{align*}
\mu_J(a\otimes[J\cap K])&\in s_*\bigl[H_*(I\cap K)\bigr],\\
\mu_L([K\cap L]\otimes b)&\in s_*\bigl[H_*(K\cap Q)\bigr].
\end{align*}
Combining the above relation, our claim follows. 
\end{example} 

\begin{rem} We note that in (2), if we use an arbitrary homology class $c\in H_*(P_{IJ})$ instead of the basic class $[I\cap J]$, then the open string product $\mu_J(c\otimes a)\in H_*(P_{IK})$ for $a\in H_*(P_{JK})$ may not be a constant path homology class coming from $H_*(I\cap K)$. 
\end{rem}

The above results may look rather technical, so we give a geometric description of the above result in terms of cycles. From this geometric point of view, the above results should be transparent. Let $\xi$ be a cycle representing a homology class $a=[\xi]\in H_*(P_{JK})$. Let $p_t: P_{JK} \rightarrow M$ be the evaluation map at $t$. For an oriented closed submanifold $I$, assume that $I$ and $p_t(\xi)$ are transversal for all $t\in[0,1]$ and let $\xi_t=\{\gamma\in \xi\mid \gamma(t)\in I\}=\xi\cap p_t^{-1}(I)$. For each path $\gamma\in \xi_t$, we split it at time $t$ into two paths and let 
\begin{equation*}
\xi_{t,[0,t]}=\{\gamma|_{[0,t]}\mid \gamma\in \xi_t\},\qquad
\xi_{t,[t,1]}=\{\gamma|_{[t,1]}\mid\gamma\in \xi_t\}. 
\end{equation*}
These are cycles in $P_{JI}$ and $P_{IK}$, respectively. 
We consider a diagram 
\begin{equation*}
P_{JK} \xleftarrow{\iota_t} P_{JI}\underset{I}{\overset{t}{\times}}
P_{IK} \xrightarrow{j_t}  P_{JI}\times P_{IK}.
\end{equation*} 
Here we have $(\iota_t)_!([\xi])=[\xi_t]$, and the map $j_t$ is given by $j_t(\gamma)=(\gamma_{[0,t]},\gamma_{[t,1]})\in \xi_{t,[0,t]}\times \xi_{t,[t,1]}$ for $\gamma\in \xi_t$ so that $\varphi_J([\xi])=(j_t)_*(\iota_t)_!([\xi])=[j_t(\xi_t)]$. 
We note that the cycle $\xi_{0,[0,1]}=\xi_0$ in $P_{I\cap J,K}$ represents $\mu_J([I\cap J]\otimes[\xi])$ in $H_*(P_{IK})$ up to sign, and the cycle $\xi_{1,[1,1]}$ is a cycle of constant paths in $I\cap K\in P_{IK}$. Since cycles $\xi_0$, $\xi_t$, and $\xi_1$ are homologous to each other in $P_{JK}$, we have that the cycles $\xi_{0,[0,1]}$, $\xi_{t,[t,1]}$, and $\xi_{1,[1,1]}$ are also homologous to each other in $P_{IK}$. This shows the first parts of (2) and (3). Similarly for second parts of (2) and (3). For (1), we observe that the cycle $\xi_{t,[0,t]}$ in $P_{JI}$ is homologous to a cycle $\xi_{0,[0,0]}$ of constant paths in $I\cap J$, and the cycle $\xi_{t,[t,1]}$ in $P_{IK}$ is homologous to a cycle $\xi_{1,[1,1]}$ of constant paths in $I\cap K$. 

In the above discussion, we observe that $[\xi_0]=[\xi_1]$ as elements in $H_*(P_{JK})$, since they are homologous through $\{\xi_t\}_{0\le t\le1}$. This means that the following compositions are equal: 
\begin{align*}
H_*(P_{JK})& \xrightarrow{(\iota_J)_!} H_*(P_{I\cap J,K}) 
\xrightarrow{(\iota_J)_*} H_*(P_{JK}), \\
H_*(P_{JK})& \xrightarrow{(\iota_K')_!} H_*(P_{I,I\cap K}) 
\xrightarrow{(\iota_K')_*} H_*(P_{JK}),
\end{align*}
where $(\iota_J)_*(\iota_J)_!([\xi])=[\xi_0]$ and $(\iota_K')_*(\iota_K')_!([\xi])=[\xi_1]$.

\begin{figure}
\begin{center}
\begin{tikzpicture}
\path (0,0) coordinate (1) node[below] {$J$};
\path (1) node[above] {$0$};
\fill (1) circle (1pt);
\path (2,0) coordinate (2) node[below] {$K$}; 
\path (2) node[above] {$1$};
\fill (2) circle (1pt);
\draw (1) -- (2);
\draw[->, >=stealth] (1,0) -- ++(0.1,0);
\draw  (1.5, 0.4) ellipse (0.3cm and 0.4cm);
\draw[->, >=stealth] (1.5,0.8) -- ++(0.1,0);
\path (1.5,0) coordinate (3) node[above] {$t$};
\path (1,0) node[below] {$\eta$};
\path (1.5, 0.8) node[above] {$\gamma$};

\draw node[text width=11cm] at (1,-1.7)
{\textsc{Figure 9.} A diagram for the $\mathbb H_*(LM)$-module structure on $H_*(P_{JK})$ parametrized by $t\in[0,1]$: 
$\bullet_t:H_*(LM)\otimes H_*(P_{JK}) \to H_*(P_{JK})$. By letting $t$ move in the interval $0\le t\le 1$, the diagram shows that left and right $\mathbb H_*(LM)$-module structures on $H_*(P_{JK})$ are the same. };

\end{tikzpicture}
\end{center}
\end{figure}

In fact, this is a special case of a $\mathbb H_*(LM)$-bimodule structure on $H_*(P_{JK})$ described by the following diagram: 
\begin{equation*}
\begin{CD}
LM\times P_{JK} @<{j_t}<< LM\underset{M}{\overset{t}{\times}}P_{JK} @>{\iota_t}>> P_{JK} \\
@V{p\times p_t}VV  @VVV @. \\
M\times M @<{\phi}<< M, @. 
\end{CD}
\end{equation*}
where $LM\underset{M}{\overset{t}{\times}}P_{JK}=(p\times p_t)^{-1}\bigl(\phi(M)\bigr)$ consists of those pairs $(\gamma,\eta)\in LM\times P_{JK}$ such that $\gamma(0)=\eta(t)$, and the map $\iota_t$ is given by $\iota_t(\gamma,\eta)=\eta_{[0,t]}\cdot\gamma\cdot\eta_{[t,1]}$. See Figure 9. For $a\in \mathbb H_*(LM)$ and $b\in H_*(P_{JK})$, the left and right $\mathbb H_*(LM)$ module structure on $H_*(P_{JK})$ is given by 
\begin{align*}
a\bullet b&=(-1)^{d|a|}(\iota_0)_*(j_0)_!,\\
b\bullet a&=(-1)^{|b|(|a|+d)}(\iota_1)_*(j_1)_!. 
\end{align*}
In particular, left and right module structures are identical, due to homotopy invariance of transfer maps. Let $[I]\in\mathbb H_*(M)\subset\mathbb H_*(LM)$ be the fundamental class of the submanifold $I$. Then with respect to the above $\mathbb H_*(LM)$ module structure,  we have 
\begin{equation}\label{LM-module structure}
\begin{aligned}
(\iota_J)_*(\iota_J)_!(b)&=\pm[I]\bullet b, \\
(\iota_K')_*(\iota_K')_!(b)&=\pm b\bullet[I]. 
\end{aligned}
\end{equation}

\subsection{Iterated open string coproducts} 

In the previous subsection, we showed that the image of the open string coproduct consists of homology classes of constant open strings. By repeatedly applying this fact, we see that the image of an iterated coproduct consists of homology classes of constant open strings living in certain submanifolds of $M$. We can systematically identify these submanifolds in which images of iterated coproducts live. 

Let $I, J, K_1,K_2,\dotsc,K_r,\dotsc$ be closed oriented submanifolds
of $M$. These are labels of vertices. An iterated coproduct is
given as a string operation associated to a sewing diagram which is a path graph with vertices $u, v_1,\dotsc, v_r,w$ directed from $u$ to $w$, with labels
$I, K_1,\dotsc,K_r, J$ such that all degree $2$ vertices are of type
(II) so that open strings are split at these vertices. Regarding this
path as a unit interval, let $0\le t_1\le t_2\le\dotsb\le t_r\le 1$ be
locations of vertices $v_1,v_2,\dotsc,v_r$. See Figure 10 (a).

We allow vertices to coincide as long as the ordering among them is preserved. Let $\vec{t}=(t_1,t_2,\dotsc,t_r)$. We consider the following diagram.
\begin{equation*}
\begin{CD}
P_{IJ} @<{\iota_{\vec{t}}}<< 
p_{\vec{t}}^{-1}(K_1\times \dotsb\times K_r) 
@>{j_{\vec{t}}}>> 
P_{IK_1}\times P_{K_1K_2}\times \dotsb\times P_{K_rJ}  \\
@V{p_{\vec{t}}}VV  @V{p_{\vec{t}}}VV  @. \\
M\times \dotsb\times M @<{\iota_1\times\dotsb\times\iota_r}<<
K_1\times \dotsb\times K_r @. 
\end{CD}
\end{equation*}
Here $p_{\vec{t}}=(p_{t_1},p_{t_2},\dotsc,p_{t_r})$ is a multiple evaluation map at time $\vec{t}=(t_1,t_2,\dots,t_r)$, and for any element $\gamma\in p_{\vec{t}}^{-1}(K_1\times \dotsb\times K_r)$, the map $j_{\vec{t}}$  is given by $j_{\vec{t}}(\gamma)=(\gamma_{[0,t_1]},\gamma_{[t_1,t_2]}, \dotsc,\gamma_{[t_r,1]})$. 
Then the iterated coproduct map is given by 
\begin{equation*}
\varphi_{K_1,\dotsc,K_r}=(j_{\vec{t}})_*\circ(\iota_{\vec{t}})_!: 
H_*(P_{IJ}) \rightarrow H_*(P_{IK_1})\otimes H_*(P_{K_1K_2})\dotsb\otimes H_*(P_{K_rJ}).
\end{equation*}
The inclusion map $\iota_1\times \cdots\times\iota_r$ has codimension $\sum_{i=1}^r(d-k_i)$, where $d=\dim M$ and $k_i=\dim K_i$. So the iterated coproduct map $\varphi_{K_1\cdots K_r}$ decreases the homological degree by the same amount. 

\begin{figure}
\begin{center}
\begin{tikzpicture}
\draw[thick] (0,0) -- (5,0);
\fill (0,0) circle (2pt);
\fill (5,0) circle (2pt);
\path (0,0) node[above] {$0$};
\path (0,0) node[below] {$u$};
\path (0,-0.6) node {$I$};
\path (5,0) node[above] {$1$};
\path (5,0) node[below] {$w$};
\path (5,-0.6) node {$J$};
\fill (1,0) circle (2pt);
\path (1,0) node[above] {$t_1$};
\path (1,0) node[below] {$v_1$};
\path (1,-0.6) node {$K_1$};
\fill (2,0) circle (2pt);
\path (2,0) node[above] {$t_2$};
\path (2,0) node[below] {$v_2$};
\path (2,-0.6) node {$K_2$};
\fill (4,0) circle (2pt);
\path (4,0) node[above] {$t_r$};
\path (4,0) node[below] {$v_r$};
\path (4,-0.6) node {$K_r$};
\draw[dashed] (2.5, 0.2) -- (3.5,0.2);
\draw[dashed] (2.5,-0.2) -- (3.5,-0.2);
\draw[->,>=stealth] (0.5,0) -- ++(0.1,0);
\draw[->, >=stealth] (1.5,0) -- ++(0.1,0);
\draw[->,>=stealth] (4.5,0) -- ++(0.1,0);

\path (-1,0) node {$(a)$};
\end{tikzpicture}
\end{center}
 
\begin{center}
\begin{tikzpicture}
 
\fill (7.5,0) circle (2pt);
\path (7,0) node[above] {$0=t_1=\cdots=t_r$};
\path (7,0) node[below] {$I\cap K_1\cap\cdots\cap K_r$};

\fill (9,0) circle (2pt);
\path (9,0) node[above] {$1$};
\path (9,0) node[below] {$J$};

\fill (11.5,0) circle (2pt);
\path (11.5,0) node[above] {$0$};
\path (11.5,0) node[below] {$I$};

\fill (13,0) circle (2pt);
\path (13.5,0) node[above] {$t_1=\cdots=t_r=1$};
\path (13.5,0) node[below] {$K_1\cap \cdots \cap K_r\cap J$};

\draw[thick] (7.5,0) -- (9,0)   (11.5,0) -- (13,0);
\draw[->,>=stealth] (8.25,0) -- ++(0.1,0);
\draw[->, >=stealth] (12.25,0) -- ++(0.1,0);

\path (4.6,0) node {$(b)$};
\path (10.8,0) node {$(c)$};

\draw node[text width=11cm] at (10,-2) 
{\textsc{Figure 10.} Deformation of sewing diagrams for iterated open string coproducts. In (b), vertices $v_1,\dots, v_r$ coincide with $u$. In (c), vertices $v_1,\dots,v_r$ coincide with $w$. The associated string operations remain the same through this deformation process. This observation provides information on the value of the iterated open string coproduct. See Theorem 4.5.};

\end{tikzpicture}
\end{center}
\end{figure}

By homotopy invariance of transfer maps, the iterated coproduct $\varphi_{K_1\dotsb K_r}$ is independent of the location of type (II) vertices $v_1,\dotsc, v_r$, that is, it is independent of the choice of parameter $\vec{t}=(t_1,t_2,\dotsc,t_r)$ with $0\le t_1\le t_2\le\cdots\le t_r\le 1$. By specializing $\vec{t}$ and using the homotopy invariance of the open string coproduct, we can obtain topological information on iterated coproducts. 

For a submanifold $N\subset K\cap L$, let $s_*:H_*(N) \rightarrow H_*(P_{KL})$ be the homomorphism induced by the inclusion map. 

\begin{thm} \label{iterated coproducts} 
The image of the iterated coproduct map 
\begin{equation*}
\varphi_{K_1\dotsb K_r}: H_*(P_{IJ}) \rightarrow 
H_*(P_{IK_1})\otimes H_*(P_{K_1K_2})\otimes\dotsb\otimes H_*(P_{K_rJ})
\end{equation*}
is contained in the subgroup 
\begin{multline*}
s_*\bigl[H_*(I\cap K_1\cap\dotsb\cap K_r)\bigr]\otimes \bigotimes^{r-1}
\bigl\{s_*\bigl[H_*(I\cap K_1\cap \dotsb\cap K_r)\bigr]\cap s_*\bigl[H_*(K_1\cap\dotsb\cap K_r\cap J)\bigr]\bigr\}\\
\otimes s_*\bigl[H_*(K_1\cap\dotsb \cap K_r\cap J)\bigr]
\end{multline*}
consisting of homology classes of constant open strings, where $s_*$'s are homomorphisms induced by suitable inclusion maps, and intersections of submanifolds are taken as transversal intersections. 
\end{thm}
\begin{proof} First we specialize $\vec{t}$ to $\vec{t}_0=(0,0,\dotsc,0)$. See the diagram in Figure 10 (b). In this case vertices $v_1,v_2,\dotsc,v_r$ coincide with $u$ and the new vertex carries the label $I\cap K_1\cap \dotsb \cap K_r$, and the space of interaction configurations is given by  $p_{\vec{t}_0}^{-1}(K_1\times\dotsb\times K_r)=P_{I\cap K_1\cap\dotsb\cap K_r,J}$. The diagram associated to this deformed case is given as follows. 
\begin{equation*}
\begin{CD}
P_{IJ} @<{\iota_0=\iota_{\vec{t}_0}}<< P_{I\cap K_1\cap\dotsb\cap K_r,J} @>{j_0=j_{\vec{t}_0}}>> 
(I\cap K_1\cap\dotsb\cap K_r)^r\times 
 P_{I\cap K_1\cap\dotsb\cap K_r,J} \\
@V{p_0}VV  @V{p_0}VV  @. \\
I @<{\iota}<<  I\cap K_1\cap\dotsb\cap K_r @. 
\end{CD}
\end{equation*}
Here the inclusion map $\iota$ has codimension $\sum_{i=1}^r(d-k_i)$. For $a\in H_*(P_{IJ})$ the effect of $(\iota_0)_!$ is given by  
\begin{equation*}
(\iota_0)_!(a)=\pm\mu_I([I\cap K_1\cap\dotsb\cap K_r]\otimes a), 
\end{equation*}
and $j_0=(p_0,p_0,\dotsc,p_0,1)\circ\phi_{r+1}$, where $\phi_{r+1}$ is the iterated diagonal map. The iterated coproduct map $\varphi$ factors as follows:
\begin{equation*}
\begin{CD}
H_*(P_{IJ})  @>{(j_0)_*(\iota_0)_!}>> H_*(I\cap K_1\cap\dotsb\cap K_r)^{\otimes r}\otimes H_*(P_{I\cap K_1\cap\dotsb\cap K_r,J})\\
@| @V{s_*\otimes\dotsb\otimes s_*\otimes (\iota_{K_r})_*}VV  \\
H_*(P_{IJ}) @>{\varphi}>> H_*(P_{IK_1})\otimes \dotsb\otimes
H_*(P_{K_{r-1}K_r})\otimes H_*(P_{K_rJ})
\end{CD}
\end{equation*}
Next, we specialize $\vec{t}$ to $\vec{t}_1=(1,1,\dotsc,1)$. See the diagram in Figure 10 (c). In this case, vertices $v_1,v_2,\dotsc,v_r$ coincide with the vertex $w$ and the new vertex carries the label $K_1\cap \dotsb\cap K_r\cap J$. The space of interaction configurations $p_{\vec{t}_1}^{-1}(K_1\times \dotsb\times K_r)$ is now $P_{I,K_1\cap\dotsb\cap K_r\cap J}$. The string operation for this degenerate sewing diagram is computed using the following diagram:
\begin{equation*}
\begin{CD}
P_{Ij} @<{\iota_1=\iota_{\vec{t}_1}}<< 
P_{I,K_1\cap\dotsb\cap K_r\cap J} @>{j_1=j_{\vec{t}_1}}>> 
P_{I,K_1\cap\dotsb\cap K_r\cap J}\times 
(K_1\cap\dotsb\cap K_r\cap J)^r \\
@V{p_1}VV  @V{p_1}VV @. \\
J @<{\iota}<<  K_1\cap \dotsb\cap K_r\cap J @.
\end{CD}
\end{equation*}
The inclusion map $\iota$ again has codimension $\sum_{i=1}^r(d-k_i)$. 
As before, for $a\in H_*(P_{IJ})$ we have $(\iota_1)_!(a)=\pm\mu_J(a\otimes [K_1\cap\dotsb\cap K_r\cap J])$, and the map $j_1$ is given by $j_1=(1,p_1,\dotsc,p_1)\circ\phi_{r+1}$. Now the iterated coproduct map $\varphi$ factors as follows: 
\begin{equation*}
\begin{CD}
H_*(P_{IJ}) @>{(j_1)_*(\iota_1)_!}>>  H_*(P_{I,K_1\cap\dotsb\cap K_r\cap J})\otimes H_*(K_1\cap\dotsb\cap K_r\cap J)^{\otimes r} \\
@|  @V{(\iota_{K_1})_*\otimes s_*\otimes\dotsb\otimes s_*}VV \\
H_*(P_{IJ}) @>{\varphi}>> 
H_*(P_{IK_1})\otimes H_*(P_{K_1K_2})\otimes \dotsb\otimes 
H_*(P_{K_rJ})
\end{CD}
\end{equation*}
Combining the above two factorizations of the iterated coproduct, we complete the proof of Theorem \ref{iterated coproducts}. 
\end{proof} 

For example, the twice iterated coproduct map 
\begin{equation*}
\varphi_{K,L}: H_*(P_{IJ}) \rightarrow H_*(P_{IK})\otimes H_*(P_{KL})
\otimes H_*(P_{LJ})
\end{equation*}
has image contained in the following tensor product:
\begin{equation*}
s_*\bigl[H_*(I\cap K\cap L)\bigr]\otimes 
\bigl(s_*\bigl[H_*(I\cap K\cap L)\bigr]\cap 
s_*\bigl[H_*(K\cap L\cap J)\bigr]\bigr)
\otimes s_*\bigl[H_*(K\cap L\cap J)\bigr]
\end{equation*}

Iterated coproduct maps are string operations associated to open
cobordisms homeomorphic to discs with one incoming open string and
many outgoing open strings. We next consider an open cobordism
$\Sigma$ with no incoming open string and only outgoing open strings
with labels $I, K_1,\dotsc,K_r,I$, where $I$ is the outermost
label. The associated string operation can be computed by evaluating
an iterated coproduct $\varphi_{K_1\dotsb K_r}$ on the fundamental
class $[I]\in H_*(P_{II})$. By using either of the factorization in
the above proof with $I=J$, we see that $\varphi_{K_1\dotsb
K_r}([I])$ is equal to
\begin{multline*}
\pm (s_*\otimes\dotsb\otimes s_*\otimes (\iota_{K_r})_*) \circ
(p_{0*}\otimes \dotsb\otimes p_{0*}\otimes 1)\circ (\phi_{r+1})_*
\bigl([I\cap K_1\cap\dotsb\cap K_r]\bigr) \\
=\pm (s_*\otimes\dotsb\otimes s_*)\circ(\phi_{r+1})_*\bigl([I\cap K_1\cap\dotsb\cap K_r]\bigr).
\end{multline*}
Thus we obtain the following proposition. 

\begin{prop}\label{only outgoing open strings} The string operation $\mu_{\Sigma}$ associated to an open cobordism $\Sigma$ homeomorphic to a disc with only outgoing open strings with labels $I,K_1,\dotsc,K_r,I$, and with values in 
\begin{equation*}
H_*(P_{IK_1})\otimes H_*(P_{K_1K_2})\otimes\dotsb\otimes H_*(P_{K_{r-1}K_r})\otimes H_*(P_{K_r,I}),
\end{equation*}
is given by the image of the diagonal map of the homology class of transversal intersection  $[I\cap K_1\cap\dotsb\cap K_r]$. Namely, 
\begin{equation*}
\mu_{\Sigma}(1)=\varphi_{K_1\dotsb K_r}([I])
=\pm (s_*\otimes \dotsb\otimes s_*)\circ(\phi_{r+1})_*\bigl([I\cap K_1\cap\dotsb\cap K_r]\bigr).
\end{equation*}
\end{prop}

\section{String operations associated to disc cobordisms}

Let $\Sigma$ be an open cobordism homeomorphic to a disc, with  incoming open strings and outgoing open strings along its boundary circle. By applying open string products to combine adjacent incoming strings and by applying open string coproducts to combine adjacent outgoing strings, the computation of the string operation $\mu_{\Sigma}$ is reduced to the computation of the string operation associated to an open cobordism homeomorphic to a disc such that incoming open strings and outgoing open strings alternate along the boundary. This is the canonical decomposition discussed in \S 2.2. 

In the previous section, we discussed properties of iterated open coproducts. In this section, we examine properties of string operations associated to disc cobordisms in which  incoming open strings and outgoing open strings alternate along their boundaries. Suppose such a disc cobordism $\Sigma$ has $m$ incoming open strings and $m$ outgoing open strings. We assume that incoming open strings are oriented clockwise and outgoing open strings oriented counter-clockwise along the boundary of $\Sigma$. Let the arcs between open strings be labeled by closed oriented submanifolds $K_1,K_2,\dotsc,K_{2m}$ in clockwise fashion so that incoming strings have odd initial labels and even terminal labels. For simplicity, we simply use indices of submanifolds $1,2,\dotsc,2m$ as labels. We denote the open string configuration space $P_{K_iK_j}$ of paths from points in $K_i$ to points in $K_j$ simply by $P_{ij}$. Thus the configuration space of incoming open strings is given by $P_{12}\times P_{34}\times\dotsb\times P_{2m-1,2m}$, and the configuration space of outgoing open strings is given by $P_{32}\times P_{54}\times\dotsb\times P_{1,2m}$. See Figure 2 in Introduction. 

The string operation associated to $\Sigma$ can be computed by decomposing it into a sequence of half-pair-of-pants, and then composing corresponding open string products or coproducts according to the decomposition. We showed in Theorem \ref{independence of decomposition} that the associated string operation $\mu_{\Sigma}$ is independent of such ordered decompositions, up to sign. Furthermore, in \S 3.2, we showed that we can compute string operations using saddle interactions where open strings interact at their internal points or at their end points, and these interactions can simultaneously take place at the same interaction point, which can be internal or end points of open strings. Using these observations, we compute the string operation $\mu_{\Sigma}$ as follows. 

Let $\vec{t}=(t_1,t_2,\dotsc,t_m)\in [0,1]^m$ be an arbitrary point in the unit $m$-dimensional cube. Consider the following diagram of simultaneous saddle interaction of $m$ open strings at the same point.
\begin{equation}\label{simultaneous interaction}
\begin{CD}
P_{12}\times P_{34}\times\dotsb\times P_{2m-1,2m} @<{\iota_{\vec{t}}}<< p_{\vec{t}}^{-1}\bigl(\phi_r(M)\bigr) 
@>{j_{\vec{t}}}>> 
P_{32}\times P_{54}\times\dotsb\times P_{1,2m} \\
@V{p_{\vec{t}}}VV   @V{p_{\vec{t}}}VV  @. \\
M\times M\times\dotsb\times M @<{\phi_m}<<  M @. 
\end{CD}
\end{equation}
Here $p_{\vec{t}}=(p_{t_1},p_{t_2},\dotsc,p_{t_m})$ is an evaluation map at specified time $\vec{t}$, and $\iota_{\vec{t}}$ is an inclusion map into the configuration space of incoming open strings from the interaction configuration subspace $p_{\vec{t}}^{-1}\bigl(\phi_r(M)\bigr)$ consisting of $m$-tuples of open strings $(\gamma_1,\gamma_2,\dotsc,\gamma_m)$ intersecting at a single point given by  $\gamma_1(t_1)=\gamma_2(t_2)=\dotsb=\gamma_m(t_m)$. The map $j_{\vec{t}}$ splits the $k$th incoming open string $\gamma_k\in P_{2k-1,2k}$ at time $t_k$ for all $k$, and recombines them so that the $[0, t_{k+1}]$ segment of $\gamma_{k+1}$ is joined by the $[t_k,1]$ segment  of $\gamma_k$ in this order to form an outgoing open string $\eta_k$. Thus, 
\begin{equation}\label{map j_t}
\begin{gathered} 
j_{\vec{t}}(\gamma_1,\gamma_2,\dotsc,\gamma_m)=(\eta_1,\eta_2,\dotsc,\eta_m), \\
\text{ \ where \ }\eta_k=(\gamma_{k+1}|_{[0,t_{k+1}]})\cdot (\gamma_{k}|_{[t_k,1]}) \text{ \ for \ }1\le k\le m,
\end{gathered}
\end{equation}
where $\gamma_{m+1}=\gamma_1$. See Figure 3 in Introduction.
The string operation $\mu_{\Sigma}$ can be computed up to sign as 
\begin{equation}\label{parametrized string operation}
\mu_{\Sigma}(\vec{t})=(j_{\vec{t}})_*\circ(\iota_{\vec{t}})_!:
H_*(P_{12}\times\dotsb\times P_{2m-1,2m}) \longrightarrow 
H_*(P_{32}\times\dotsb\times P_{1,2m}).
\end{equation}
By the homotopy invariance of transfer maps, the above homomorphism is independent of the choice of $\vec{t}\in [0,1]^m$. 

By specializing $\vec{t}$ to certain boundary points of $[0,1]^m$, we can expose various topological properties of the above string operation $\mu_{\Sigma}$. In particular, we show that the image of $\mu_{\Sigma}$ consists of homology classes of constant open strings lying in certain submanifolds which we now precisely determine.  
To describe the result, let $E=\{\varepsilon=(\varepsilon_1,\varepsilon_2,\dotsc,\varepsilon_m) \mid \varepsilon_k=\pm 1, \text{ for }1\le k\le m\}$ be the totality of sequences of $\pm1$ of length $m$. This set $E$ is the set of coordinates of corners of $m$-dimensional cube $[0,1]^m$. For $\varepsilon\in E$, let 
\begin{equation*}
K_{\varepsilon}=\bigcap_{j=1}^{m}K_{2j-1+\varepsilon_j}.
\end{equation*}
For each $1\le k\le m$, 
let $E_k=\{\varepsilon\in E \mid \varepsilon_k=1, \varepsilon_{k+1}=0\}$, where for $k=m$ we set $\varepsilon_{m+1}=\varepsilon_1$. This is the set of corner points of $[0,1]^m$ which affect the $H_*(P_{2k+1,2k})$ factor in the image of the string operation $\mu_{\Sigma}$, as we will show. For $\varepsilon\in E_k$, 
\begin{equation*}
K_{\varepsilon}=\bigl[K_{2k}\cap K_{2k+1}\cap\bigl(\!\!\!\bigcap_{\substack{1\le j\le m \\j\not=k,k+1}}\!\!\!K_{2j-1+\varepsilon_j}\bigr)\bigr]\subset K_{2k}\cap K_{2k+1}
\subset P_{2k+1,2k}.
\end{equation*}
Here $P_{2m+1,2m}=P_{1,2m}$ when $k=m$. For $\varepsilon\in E_k$, let $s_{\varepsilon}:K_{\varepsilon} \rightarrow P_{2k+1,2k}$ be the inclusion map. We define a subgroup $S_{2k+1,2k}\subset H_*(P_{2k+1,2k})$ by 
\begin{equation}\label{S group}
S_{2k+1,2k}=\!\bigcap_{\varepsilon\in E_k} \!(s_{\varepsilon})_*\bigl[H_*(K_{\varepsilon})\bigr] \subset H_*(P_{2k+1,2k}).
\end{equation}

\begin{thm} \label{disc operation} Let $\Sigma$ be an open cobordism homeomorphic to a disc with alternating $m$ incoming and $m$ outgoing open strings whose free arcs are labeled by $K_1,K_2,\dotsc,K_{2m}$ as above. Then the associated string operation
\begin{equation*}
\mu_{\Sigma}: \bigotimes_{k=1}^mH_*(P_{2k-1,2k}) 
\longrightarrow \bigotimes_{k=1}^m H_*(P_{2k+1,2k})
\end{equation*}
has its image in the subgroup of homology classes of constant open strings given by 
\begin{equation*}
\textup{Im}\,(\mu_{\Sigma})\subset 
\bigotimes_{k=1}^m S_{2k+1,2k}.
\end{equation*}
\end{thm}
\begin{proof} For $\varepsilon\in E\subset [0,1]^m$, in the interaction associated to the string operation $\mu_{\Sigma}(\varepsilon)$, $m$ incoming open strings meet simultaneously at one of the end points of $\gamma_j$ labeled by a submanifold with index $2j-1+\varepsilon_j$ for $1\le j\le m$. At this moment of interaction, the intersection point belongs to the submanifold $K_{\varepsilon}$. Right after the moment of interaction, this configuration splits into $m$ outgoing open strings $\eta_1,\eta_2,\dotsc,\eta_m$ by $j_{\varepsilon}$. When $\varepsilon\in E_k\subset E$, the diagram relevant to this interaction is the following one:
\begin{equation}\label{e-interaction}
\!\!\!\!\!\!\!\!\!\!\!\!\!\!\!\!\!\!\!\!\!\!\!\!
\begin{CD}
P_{12}\times\dotsb\times P_{2k-1,2k}\times P_{2k+1,2k+2}\times\dotsb\times P_{2m-1,2m} @<{\iota_{\varepsilon}}<< 
p_{\varepsilon}^{-1}\bigl(\phi_m(K_{\varepsilon})\bigr) @>{j_{\varepsilon}}>> 
P_{32}\times\dotsb\times P_{1,2m} \\
@V{p_{\varepsilon}}VV  @V{p_{\varepsilon}}VV  @. \\
K_{1+\varepsilon_1}\times \dotsb\times K_{2k}\times K_{2k+1}\times\dotsb\times K_{2m-1+\varepsilon_m} @<{\phi_m}<< 
K_{\varepsilon} @.  \\
@VVV @VVV @. \\
M\times\dotsb\times M\times M\times \dotsb\times M @<{\phi_m}<< 
M @. 
\end{CD}
\end{equation}
Here the bottom square is a pull-back diagram which appears when type (2) deformation further degenerates into type (3) deformation discussed in \S 3.2. String operation remains the same throughout this deformation process.  

When we examine the map $j_{\vec{t}}$ in \eqref{map j_t} for $\vec{t}=\varepsilon$, we see that for $\varepsilon\in E_k$ the $k$\,th outgoing open string $\eta_k=(\gamma_{k+1}|_{[0,0]})\cdot(\gamma_{k}|_{[1,1]})$ is a constant open string lying in $K_{\varepsilon}$. See Figure 4 in Introduction. Thus $j_{\varepsilon}$ for $\varepsilon\in E_k$ factors as follows: 
\begin{equation*}
\xymatrix{
p_{\varepsilon}^{-1}\bigl(\phi_m(K_{\varepsilon})\bigr) \ar[r]^{j_{\varepsilon}\ \ \ \ \ \ \ \ \ \ \ \ \ \ } \ar[dr] & P_{32}\times\dotsb\times P_{2k+1,2k}\times\dotsb\times P_{1,2m} \\
& P_{32}\times\dotsb\times K_{\varepsilon}\times\dotsb\times P_{1,2m}
\ar[u]
}
\end{equation*}
Hence the $k$-th factor of the image of the string operation $\mu_{\Sigma}(\varepsilon)$ in the homology group $\bigotimes_{k=1}^mH_*(P_{2k+1,2k})$ comes from $H_*(K_{\varepsilon})$ for any $\varepsilon\in E_k$. Hence the $k$-th factor of the image of $\mu_{\Sigma}(\varepsilon)$ is actually contained in $S_{2k+1,2k}$ given in the intersection \eqref{S group}. Since we are taking the tensor product over a field, applying the above argument for each tensor factor, we see that the image of the string operation is contained in $\bigotimes_{k=1}^mS_{2k+1,2k}$, which consists of homology classes of constant paths. This completes the proof. 
\end{proof} 

We describe a simple case in which the string operation $\mu_{\Sigma}$ vanishes. 

\begin{cor}\label{vanishing of string operation}  With the above notation, if $K_{\varepsilon}=\emptyset$ for some $\varepsilon\in E$, then the associated string operation $\mu_{\Sigma}$ vanishes. 
\end{cor}
\begin{proof} If $\varepsilon\in E$ is such that $\varepsilon\ne(0,0,\dotsc,0), (1,1,\dotsc,1)$, then there exists $k$ such that $\varepsilon_k=1$ and $\varepsilon_{k+1}=0$, where we set $\varepsilon_{m+1}=\varepsilon_1$ to take care of the case $(0,\varepsilon_2,\dotsc,\varepsilon_{m-1},1)$. For this $k$, we have $S_{2k+1,2k}\subset s_*\bigl[H_*(K_{\varepsilon})\bigr]=0$ since $K_{\varepsilon}=\emptyset$ by hypothesis. Hence $\mu_{\Sigma}\equiv0$.

When $\varepsilon=(0,0,\dots,0)$, in the corresponding simultaneous saddle interaction, all incoming open strings meet at their tail vertices, and this interaction point lies in $K_{(0,\dots,0)}=\bigcap_{j=1}^m K_{2j-1}$. See Figure 4 in Introduction. If $K_{(0,\dots,0)}=\emptyset$, then in the interaction diagram \eqref{e-interaction} the space of interaction configurations $p_{(0,\dots,0)}^{-1}\bigl(\phi_m(K_{(0,\dots,0)})\bigr)$ is an empty set. Hence the associated string operation $\mu_{\Sigma}$ vanishes. 

Similarly for the case $\varepsilon=(1,1,\dots,1)$ in which all incoming open strings meet at their head vertices in the submanifold $K_{(1,\dots,1)}$. 
This completes the proof. 
\end{proof} 

We apply Theorem \ref{disc operation} in the simplest nontrivial case $m=2$, which describes the saddle operation introduced in \cite{T5}. 

\begin{cor}\label{m=2 example} Let $\Sigma$ be an open cobordism homeomorphic to a disc with two incoming and two outgoing open strings arranged alternately, with free arc labels given by oriented closed submanifolds $K_1,K_2,K_3,K_4$. Then the associated string operation
\begin{equation*}
\mu_{\Sigma}: H_*(P_{12})\otimes H_*(P_{34}) \longrightarrow 
H_*(P_{32})\otimes H_*(P_{14})
\end{equation*}
has image contained in the subgroup 
\begin{equation*}
\textup{Im}\,\mu_{\Sigma}\subset 
s_*\bigl[H_*(K_3\cap K_2)\bigr]\otimes 
s_*\bigl[H_*(K_1\cap K_4)\bigr],
\end{equation*}
where $s$ is the inclusion map from the subspace of the constant paths. 
\end{cor} 

The next case is for $m=3$. 

\begin{cor} \label{m=3 example} Let $\Sigma$ be a disc cobordism with $3$ incoming and $3$ outgoing open strings arranged alternately along its boundary with free arc labels $K_1,K_2,\dotsc,K_6$. Then the associated string operation 
\begin{equation*}
\mu_{\Sigma}: H_*(P_{12})\otimes H_*(P_{34})\otimes H_*(P_{56})
\longrightarrow 
H_*(P_{32})\otimes H_*(P_{54})\otimes H_*(P_{16})
\end{equation*}
has its image in the subgroup of homology classes of constant paths given as follows\textup{:} 
\begin{multline*}
\textup{Im}\,\mu_{\Sigma}\subset 
\bigl(s_*\bigl[H_*(K_3\cap K_2\cap K_5)\bigr] \cap 
s_*\bigl[H_*(K_3\cap K_2\cap K_6)\bigr]\bigr) \\
\otimes \bigl(s_*\bigl[H_*(K_5\cap K_4\cap K_1)\bigr] \cap s_*\bigl[H_*(K_5\cap K_4\cap K_2)\bigr]\bigr)  \\
\otimes \bigl(s_*\bigl[H_*(K_1\cap K_6\cap K_3)\bigr]\cap s_*\bigl[H_*(K_1\cap K_6\cap K_4)\bigr]\bigr),
\end{multline*}
where $s_*$'s are induced by obvious inclusion maps. 
\end{cor} 

\bibliography{bibliography}
\bibliographystyle{plain}

\end{document}